\documentclass[11pt]{article}
\usepackage{amsmath,amssymb,amsthm,bm,geometry,graphicx,booktabs,array,caption}
\usepackage{placeins}
\usepackage{url}
\graphicspath{{./}}

\newtheorem{theorem}{Theorem}
\newtheorem{proposition}[theorem]{Proposition}
\newtheorem{corollary}[theorem]{Corollary}
\newtheorem{remark}[theorem]{Remark}

\title{Alias Geometry and Density Compensation in Perturbed Fourier Sampling}
\author{
  Erli Wind-Andersen\\[-1pt]
  \small Department of Genetics and Genomic Sciences\\[-1pt]
  \small Icahn School of Medicine at Mount Sinai\\[-1pt]
  \small New York, NY, USA\\[-1pt]
  \small \texttt{erli.wind-andersen@mssm.edu}
}
\date{}

\begin{document}
\maketitle
\vspace{-2.2em}

\begin{abstract}
On a uniform grid with equal weights, Fourier modes that remain distinct on the grid are exactly orthogonal. When the nodes are perturbed, density compensation readjusts the weights. We ask how much of the resulting error reweighting can remove, measuring it by the Frobenius norm of the Gram defect, whose square is proportional to the mean-square error of density-compensated adjoint reconstruction. The answer is governed by aliasing. The grid cannot distinguish certain frequency differences, so a weight correction must treat each such alias class as a single number, while a node displacement affects every frequency difference in the class separately. To first order, reweighting removes the class average but not the spread. The irreducible error is a sum of covariance matrices, one for each class of colliding frequency differences. If no differences collide, every small node perturbation can be repaired exactly with positive weights. With collisions, a repair that works to first order can fail at second order. At the two-dimensional DFT, row and column shears are each exactly repairable, but their combination is not, and we compute the obstruction. In one dimension, oversampling protects a band of low-frequency node motions, and under random jitter optimal reweighting asymptotically removes a fraction \(1-\kappa^2(2-\kappa)^2\) of the expected squared first-order error at oversampling ratio \(\kappa\in[1,2]\). At critical sampling, node and weight perturbations are orthogonal: no reweighting helps, Voronoi weights can never reduce the error, and for the one-dimensional DFT the remaining error is a log-Vandermonde energy with the circular \(H^{1/2}\) energy as its continuum limit.
\end{abstract}
\vspace{-0.8em}


\vspace{0.3em}

\vspace{0.3em}
\noindent\textbf{Reproducibility materials.}
Code that reproduces every numerical value and figure in this paper, together with symbolic and numerical checks of its identities, is available at \url{https://github.com/circlespie/alias-geometry-fourier-sampling}.

\section{Introduction}

Nonuniform Fourier sampling comes with two geometric ingredients: where the samples are
taken and how they are weighted. When the nodes move away from a configuration at which the
resolved Fourier modes are exactly orthogonal, one can try to compensate by changing the
weights. Which part of the node-induced Gram perturbation can diagonal reweighting eliminate,
which part necessarily remains, and what feature of the sampling geometry distinguishes the
two? We study this local compensability rather than the stability of a perturbed sampler or
the optimization of weights for fixed nodes.

Perturbations of uniform Fourier grids have been studied through the stability and
convergence of interpolation and quadrature \cite{AustinTrefethen2017,YuTownsend2023},
including the conditioning of nonuniform Fourier matrices and the behavior of quadrature
weights. A complementary literature optimizes density-compensation weights for prescribed
sampling nodes, using exact moment conditions or minimization of the weighted Fourier Gram
defect \cite{SedaratNishimura2000,KircheisPottsDCF2023,KircheisPottsFrontiers2023}. This is also the usual setting in
non-Cartesian MRI, where compensation is designed for a prescribed \(k\)-space trajectory
\cite{PipeMenon1999,Rosenfeld1998}. Classical node-perturbation theory, from Kadets' \(1/4\)
theorem for exponential systems \cite{Kadec1964} to Gr\"ochenig's adaptive-weight theory of
irregular sampling \cite{Groechenig1993}, asks how far nodes may move before stability is lost;
the present theory is infinitesimal and measures in Frobenius norm what moving weights can
recover at an exact rule. In frame language, exact rules are precisely the diagonal scalings of
the Fourier frame \(\{e^{i\ell\cdot x_k}\}\) to a Parseval frame, so the question is the local
scalability of perturbed Fourier frames \cite{KutyniokOkoudjouPhilippTuley2013}.

Node motion and weight motion are coupled, but they are not interchangeable. Changing a
weight only changes the coefficient attached to a sampled Fourier character; moving a node
changes its phase, and therefore acts through the actual frequency difference between the two
modes being compared. At an exact uniform rule supported on a finite subgroup of
\(\mathbb T^d\), frequency differences that induce the same character on the sampling group
form an alias class, and a single weight correction must respond to all of them at once:
reweighting acts at the resolution of alias classes, while node motion still resolves the
individual differences inside each class (Figure~\ref{fig:alias-class-geometry}). In
character coordinates the first-order optimization separates into independent
low-dimensional least-squares problems, one for each class.

If \(C_r\) denotes the multiplicity-weighted covariance of the raw frequency differences in
class \(r\), then after the optimal first-order weight response has been removed,
\begin{equation}
\label{eq:intro-covariance}
q(\epsilon)
=
\sum_r
\widehat\epsilon_r^{\,*}C_r\widehat\epsilon_r.
\end{equation}
The rank of \(C_r\) counts the damaging polarization directions in class \(r\), while
\(\ker C_r\) consists exactly of the node polarizations that remain first-order repairable.
Two distinct raw differences \(p,p'\in\Lambda-\Lambda\) form a \emph{collision} when they
induce the same character on the sampling set, equivalently when \(p-p'\) lies in the annihilator lattice \(L_{\mathcal G}\) of the grid (Section~\ref{sec:character-classes}).
Singleton classes have \(C_r=0\); collisions create the obstruction.

If the raw difference set is injective modulo the sampling lattice, every class is a singleton
and \(q\equiv0\). More is true: sufficiently small node perturbations then admit nearby
positive normalized weights for which the rule remains exactly orthogonal, so exact
reweightability persists locally.

Collisions raise a question that first-order theory alone cannot answer. A direction in
\(\ker q\) is infinitesimally repairable, but it need not lie on any nearby family of exact
rules. When the representatives in each collision class are affinely independent, the
exactness equations form a smooth local manifold whose tangent space is precisely \(\ker q\).
Affine relations can break this equivalence: at the tensor DFT, the first-order kernel
contains simultaneous row and column shears, yet generic mixtures fail a second-order
compatibility condition, so the true tangent cone is strictly smaller.

There is a parallel here with the deformation theory of complex Hadamard
matrices. At an exact square Fourier rule, \(A\) is a complex Hadamard matrix, so that
\(A/\sqrt M\) is unitary, and the distinction between first-order admissible directions and genuine
nearby exact families is classical in that setting
\cite{Dita2004,TadejZyczkowski2006,TadejZyczkowski2008,Banica2014,BarrosSaBengtsson2013}.
Our perturbations are more structured: they arise from moving Fourier nodes, with diagonal
reweighting available simultaneously. The tensor shear example gives a Fourier-node
realization of the same formal-versus-integrable distinction, with the affine relations inside an
alias class supplying the second-order obstruction. The statement concerns the node-restricted
deformations, with phases \(\theta_{\ell k}=\ell\cdot\epsilon_k\); we make no claim about
whether the mixed directions integrate inside the full complex-Hadamard variety.

In one dimension the geometry is completely explicit. For \(m\) consecutive Fourier modes
sampled at \(s\) equally spaced nodes, the collided classes occupy a high-frequency window,
while the low-frequency node modes
\begin{equation}
\label{eq:intro-repairable-band}
|n|\le s-m
\end{equation}
are exactly the first-order repairable directions. This yields a closed-form deficit law, an
exact mean-square prediction for iid jitter, and an exact price for prescribed Voronoi weights.
At critical sampling, Voronoi weights can never reduce the first-order Frobenius defect, and
they increase it unless the perturbation is confined to the rigid mode and, for even \(s\), the Nyquist mode. In the
proportional regime \(s/m\to\kappa\in[1,2]\), the optimally repairable fraction of
the expected first-order Gram energy tends to
\begin{equation}
\label{eq:intro-jitter-limit}
1-\kappa^2(2-\kappa)^2,
\end{equation}
so oversampling enlarges the repairable sector continuously from rigid translation at critical
sampling to complete local reweightability at the difference-set-injective threshold.

At the critical square endpoint the loss of weight freedom is stronger still.
For every positive exact square Fourier rule, the node- and weight-generated
first-order Gram perturbations are Frobenius-orthogonal,
\[
P_YX=0,
\]
so diagonal reweighting cannot remove any part of the first-order node defect.
For the one-dimensional square DFT, the resulting irreducible energy is the
classical circular log-Vandermonde Hessian, whose square-band limit is the Douglas
(\(\dot H^{1/2}\)) energy on the circle \cite{Douglas1931,DiNezzaPalatucciValdinoci2012}.
The square endpoint is not merely the case of zero mean-square repair in the
one-dimensional law; it is a general first-order orthogonality phenomenon.

Our analysis concerns a local compensation problem around exact weighted
Fourier rules. We optimize the first-order defect of the frame operator in
Frobenius norm. The block decomposition is derived for uniform finite-subgroup
references, while the nonlinear results require additional hypotheses on the
exact weight and on the affine geometry of the collision classes. Spectral
objectives, already-inexact reference rules, and reconstruction-specific
criteria such as least-squares error or conditioning are separate questions.
First-order alias geometry is the object of study; the second-order condition of
Section~\ref{sec:affine-relations} is offered as the entry point for a higher-order analysis,
not as a partial attempt at one.

Sections~\ref{sec:node-vs-weight}--\ref{sec:alias-classes} set up the reduction and prove the
alias-covariance decomposition, and Sections~\ref{sec:integrability}--\ref{sec:tensor-twist}
separate first-order repairability from exact deformability. Sections~\ref{sec:one-dim}--\ref{sec:higher-dim}
work out the one-dimensional, critical, and higher-dimensional geometries explicitly.

\section{Node motion versus weight motion}
\label{sec:node-vs-weight}

\subsection{Exact weighted Fourier rules}

Let \(\Lambda\subset\mathbb Z^d\) be a finite set of \(M\) frequencies, and let
\(x_1,\ldots,x_S\in\mathbb T^d\) be the sampling nodes. The associated Fourier matrix is
\[
A_{\ell k}(x)=e^{i\ell\cdot x_k},
\qquad \ell\in\Lambda,\quad k=1,\ldots,S.
\]
We attach positive normalized weights \(v_k\), collected in
\[
W(v)=\operatorname{diag}(v_1,\ldots,v_S),
\qquad
v_k>0,
\qquad
\sum_{k=1}^S v_k=1,
\]
and follow the weighted frame operator
\[
G(x,v)=A(x)W(v)A(x)^*,
\qquad
G_{\ell q}(x,v)
=
\sum_{k=1}^S v_k e^{i(\ell-q)\cdot x_k}.
\]
The rule is \emph{exact} when \(G(x,v)=I_M\): the resolved modes are orthogonal after
weighting, equivalently the weighted node set integrates every nonzero frequency difference
\(\ell-q\) to zero. We now move the nodes away from an exact configuration and allow the
weights to move with them. Both changes affect \(G\); the question is whether the change
produced by the weights can cancel the change produced by the nodes.

\subsection{First variation}

Let \((x^0,v^0)\) be an exact rule, and perturb both parts of it:
\[
x_k(t)=x_k^0+t\epsilon_k+o(t),
\qquad
v_k(t)=v_k^0+t\zeta_k+o(t),
\qquad
\sum_k\zeta_k=0.
\]
To first order the two effects add, so \(G(t)=I+tH+o(t)\) with \(H=X\epsilon+Y\zeta\): the
term \(X\epsilon\) is created by the node motion, and \(Y\zeta\) describes everything an
admissible first-order change of weights can do in response. For a pair of resolved
frequencies \(\ell,q\in\Lambda\), write \(p=\ell-q\). Then
\begin{equation}
(Y\zeta)_{\ell q}
=
\sum_{k=1}^S \zeta_k e^{ip\cdot x_k^0},
\label{eq:weight-first-variation}
\end{equation}
whereas
\begin{equation}
(X\epsilon)_{\ell q}
=
i\sum_{k=1}^S v_k^0\,
(p\cdot\epsilon_k)\,
e^{ip\cdot x_k^0}.
\label{eq:node-first-variation}
\end{equation}
A change of weights alters only the coefficient attached to the character
\(e^{ip\cdot x}\), while a node displacement also sees the frequency difference \(p\) itself
through \(p\cdot\epsilon_k\). This difference is what will make some node errors impossible
to reproduce by reweighting alone.

\begin{proposition}[First variation]
\label{prop:first-variation}
Let \((x^0,v^0)\) be an exact positive normalized rule, perturbed as above. Then
\[
G(x(t),v(t))
=
I+t\bigl(X_{v^0}\epsilon+Y\zeta\bigr)+o(t),
\]
with \(Y\) and \(X_{v^0}\) given by \eqref{eq:weight-first-variation} and
\eqref{eq:node-first-variation}. Every zero-sum \(\zeta\) is realized by a positive normalized
weight path for sufficiently small \(|t|\).
\end{proposition}

\begin{proof}
Differentiating the defining expression for \(G_{\ell q}\) gives
\eqref{eq:weight-first-variation}--\eqref{eq:node-first-variation}. For any zero-sum
\(\zeta\), the path \(v(t)=v^0+t\zeta\) remains normalized and, by strict positivity of
\(v^0\), admissible for sufficiently small \(|t|\).
\end{proof}

\subsection{Eliminating the weight directions}

Now fix the node motion \(\epsilon\) and ask how well the weights can respond. Throughout we
regard Hermitian matrices as a real Hilbert space under
\(\langle A,B\rangle_F=\operatorname{Re}\operatorname{tr}(A^*B)\). Normalized weight variations
lie in \(\mathcal Z=\{\zeta\in\mathbb R^S:\mathbf 1^T\zeta=0\}\), so the admissible corrections
form the subspace \(Y\mathcal Z\), and the best first-order compensation projects
\(X\epsilon\) onto it. With \(P_Y\) the Frobenius-orthogonal projector onto \(Y\mathcal Z\),
the part no weight change can remove has squared size
\begin{equation}
\label{eq:q-definition}
q_{v^0}(\epsilon)
=
\bigl\|(I-P_Y)X_{v^0}\epsilon\bigr\|_F^2
=
\min_{\zeta\in\mathcal Z}
\|X_{v^0}\epsilon+Y\zeta\|_F^2 .
\end{equation}
Thus \(q_{v^0}(\epsilon)=0\) precisely when admissible reweighting cancels the first-order
node defect; otherwise \(\sqrt{q_{v^0}(\epsilon)}\) is the sharp residual norm.

The Frobenius objective has two natural interpretations. Writing \(B=A W^{1/2}\), the
frame-potential excess satisfies
\begin{equation}
\label{eq:frobenius-frame-potential}
\operatorname{FP}(B)-M
=
\|G-I\|_F^2 ,
\end{equation}
and if the resolved coefficient vector \(c\) is isotropic with
\(\mathbb E[cc^*]=\sigma_c^2 I\), then
\begin{equation}
\label{eq:frobenius-adjoint-mse}
\mathbb E\|(G-I)c\|_2^2
=
\sigma_c^2\|G-I\|_F^2 .
\end{equation}
So \(q_{v^0}(\epsilon)\) measures the irreducible first-order growth of both the
frame-potential excess \cite{BenedettoFickus2003} and the mean-square adjoint reconstruction
error after the weights have been used optimally. In particular, \(q_{v^0}(\epsilon)=0\) also
implies
\begin{equation}
\label{eq:frobenius-spectral-consequence}
\|G(t)-I\|_{\mathrm{op}}=o(t),
\qquad
\operatorname{cond}_2 G(t)=1+o(t)
\end{equation}
along an optimally compensated first-order path.

Because \(P_Y\) is linear, \(q_{v^0}\) is a quadratic form in the node displacement. Give the
node space the reference-weighted inner product
\begin{equation}
\label{eq:node-inner-product}
\langle\epsilon,\phi\rangle_{v^0}
=
\sum_{k=1}^S v_k^0\,\epsilon_k\cdot\phi_k ,
\end{equation}
for which the finite-group Fourier transform is an isometry at the uniform reference. With
\(X_{v^0}^{\sharp}\) the adjoint of the node-variation map, defined by
\(\langle X_{v^0}\epsilon,H\rangle_F=\langle\epsilon,X_{v^0}^{\sharp}H\rangle_{v^0}\),
\[
q_{v^0}(\epsilon)
=
\langle \epsilon,Q_{\rm red}\epsilon\rangle_{v^0},
\qquad
Q_{\rm red}
=
X_{v^0}^{\sharp}(I-P_Y)X_{v^0},
\]
and \(\ker Q_{\rm red}\) is exactly the space of node directions that reweighting cancels to
first order. For a general exact rule this is the complete local variable-projection reduction
\cite{GolubPereyra1973}; the Fourier structure enters when we ask what \(Q_{\rm red}\) looks
like for a regular sampling rule.

\begin{proposition}[Local reduction and fixed-node optimum]
\label{prop:local-reduction}
Let \((x^0,v^0)\) be an exact positive normalized rule, and for a prescribed node tangent
\(\epsilon\) let \(q_{v^0}(\epsilon)\) be as in \eqref{eq:q-definition}. Then
\[
\inf_{v(\cdot)}
\lim_{t\downarrow0}
\frac{
\|G(x^0+t\epsilon+o(t),v(t))-I\|_F
}{t}
=
q_{v^0}(\epsilon)^{1/2},
\]
where the infimum is over differentiable positive normalized paths
\(v(t)=v^0+t\zeta+o(t)\) with \(\zeta\in\mathcal Z\).

If \(Y|_{\mathcal Z}\) is injective, define
\[
F_*(x)
:=
\min_{\mathbf 1^Tv=1}
\|A(x)W(v)A(x)^*-I\|_F^2,
\qquad
V_*(x):=F_*(x)^{1/2},
\]
with the minimum taken over real normalized weights. Then, for \(x\) sufficiently near
\(x^0\), the minimizer \(v_*(x)\) is unique, smooth, and positive. Along
\(x(t)=x^0+t\epsilon+o(t)\), one has
\[
F_*(x(t))
=
t^2q_{v^0}(\epsilon)+o(t^2),
\qquad
V_*(x(t))
=
t\,q_{v^0}(\epsilon)^{1/2}+o(t),
\]
and
\[
Dv_*(x^0)[\epsilon]
=
-\,Y^\dagger X_{v^0}\epsilon .
\]
\end{proposition}

Here \(Y^\dagger\) denotes the Moore--Penrose pseudoinverse of \(Y|_{\mathcal Z}\).

\begin{proof}
The first claim is the orthogonal least-squares identity
\[
\min_{\zeta\in\mathcal Z}
\|X_{v^0}\epsilon+Y\zeta\|_F
=
\|(I-P_Y)X_{v^0}\epsilon\|_F ,
\]
and Proposition~\ref{prop:first-variation} shows that every \(\zeta\in\mathcal Z\) is
realized by a positive normalized path for sufficiently small \(t\).

For the second claim, write
\[
b(x)=A(x)W(v^0)A(x)^*-I .
\]
The fixed-node residual is \(b(x)+Y_x\zeta\). Injectivity of \(Y_{x^0}|_{\mathcal Z}\)
persists nearby, so the least-squares minimizer is uniquely
\[
\zeta_*(x)
=
-\,Y_x^\dagger b(x),
\qquad
v_*(x)
=
v^0+\zeta_*(x),
\]
and depends smoothly on \(x\). Since \(v^0>0\), it remains positive nearby. Finally,
\[
b(x(t))
=
tX_{v^0}\epsilon+o(t),
\qquad
P_{Y_x}\to P_Y,
\]
which gives the two asymptotic expansions; differentiating the minimizer at \(x^0\) gives the
stated formula for \(Dv_*(x^0)[\epsilon]\).
\end{proof}

\begin{corollary}[Universal first-order floor]
\label{cor:universal-first-order-floor}
Assume \(Y|_{\mathcal Z}\) is injective, and let \(x(t)=x^0+t\epsilon+o(t)\). Then for every
choice of positive normalized weights \(v(t)\), with no continuity or differentiability
assumption,
\[
\liminf_{t\downarrow0}
\frac{\|G(x(t),v(t))-I\|_F}{t}
\ge
q_{v^0}(\epsilon)^{1/2} .
\]
If \(v(t)=v^0+t\zeta+o(t)\) with \(\zeta\in\mathcal Z\), then equality in the first-order
coefficient holds if and only if
\[
Y\zeta
=
-P_YX_{v^0}\epsilon .
\]
\end{corollary}

\begin{proof}
By Proposition~\ref{prop:local-reduction}, for \(x\) sufficiently near \(x^0\), the minimum
over real normalized weights is attained by a positive weight vector and satisfies
\[
V_*(x(t))
=
t\,q_{v^0}(\epsilon)^{1/2}+o(t).
\]
Hence every positive normalized choice \(v(t)\) obeys
\[
\|G(x(t),v(t))-I\|_F
\ge
V_*(x(t)),
\]
which gives the liminf bound. For a differentiable weight path,
Proposition~\ref{prop:first-variation} gives first variation
\(X_{v^0}\epsilon+Y\zeta\). Equality with the least-squares minimum occurs exactly when its
weight-generated component cancels the projection of \(X_{v^0}\epsilon\) onto
\(Y\mathcal Z\), namely \(Y\zeta=-P_YX_{v^0}\epsilon\).
\end{proof}

\begin{remark}[Injectivity of the weight map at subgroup references]
\label{rem:weight-injectivity}
For the finite-subgroup references of Section~\ref{sec:alias-classes}, the entries of
\(Y\zeta\) are the Fourier coefficients of \(\zeta\) at the character classes met by
\(\Delta=\Lambda-\Lambda\), and \(\mathbf 1^T\zeta=0\) removes the zero class. Hence
\(Y|_{\mathcal Z}\) is injective exactly when \(\Delta\) meets every nonzero character class.
This is the uniqueness hypothesis used in Theorem~\ref{thm:exactness-manifold} and
Propositions~\mbox{\ref{prop:affine-relation-obstruction}--\ref{prop:tensor-twist}}.
\end{remark}

\section{Alias classes and the irreducible defect}
\label{sec:alias-classes}

\subsection{Character classes}\label{sec:character-classes}

Now specialize to an exact rule whose nodes form a finite subgroup
\(\mathcal G\subset\mathbb T^d\), with uniform weights \(1/S\). The first-variation
formulas \eqref{eq:weight-first-variation}--\eqref{eq:node-first-variation} already show
the relevant symmetry. If two frequency differences \(p,p'\in\mathbb Z^d\) satisfy
\[
p-p'\in L_{\mathcal G},
\qquad
L_{\mathcal G}
=
\{h\in\mathbb Z^d:
e^{ih\cdot x}=1\ \text{for every }x\in\mathcal G\},
\]
then
\[
e^{ip\cdot x}=e^{ip'\cdot x}
\qquad
\text{for every }x\in\mathcal G.
\]
Thus the sampling grid cannot distinguish \(p\) from \(p'\) at the level of the character
appearing in the weight and node variations. We group such differences into the same
\emph{alias class}, and write \([p]\) for the residue class of \(p\) modulo
\(L_{\mathcal G}\). The raw differences and their classes are
\begin{equation}
\label{eq:difference-classes}
\Delta=\Lambda-\Lambda,\qquad
\Delta_r=\{p\in\Delta:p\equiv r\pmod{L_{\mathcal G}}\}.
\end{equation}

For the uniform reference \(v_k^0=1/S\), this class structure appears directly in the
first variation. Set \(\eta_k=S\zeta_k\), and for a scalar or vector field \(f\) on
\(\mathcal G\) write
\[
\widehat f_r
=
\frac1S\sum_{k=1}^S f_k e^{ir\cdot x_k^0}.
\]
If \(p=\ell-q\), then the formulas of Section~\ref{sec:node-vs-weight} reduce to
\[
(Y\zeta)_{\ell q}
=
\widehat\eta_{[p]},
\qquad
(X\epsilon)_{\ell q}
=
i\,p\cdot\widehat\epsilon_{[p]},
\]
and hence
\begin{equation}
H_{\ell q}
=
\widehat\eta_{[p]}
+i\,p\cdot\widehat\epsilon_{[p]}.
\label{eq:classwise-first-variation}
\end{equation}
The important asymmetry is now explicit: the weight correction \(\widehat\eta_{[p]}\) depends only on
the alias class, whereas the node term still retains the particular raw difference \(p\).
Thus every \(p\) in the same class must share one weight correction, even though different
representatives generally demand different corrections.

The Frobenius norm now groups naturally by alias class. If
\[
N_p
=
\#\{(\ell,q)\in\Lambda^2:\ell-q=p\},
\]
then
\begin{equation}
\|H\|_F^2
=
\sum_r\sum_{p\in\Delta_r}
N_p\,
\bigl|
\widehat\eta_r+i\,p\cdot\widehat\epsilon_r
\bigr|^2.
\label{eq:classwise-frobenius}
\end{equation}
The minimization separates into independent scalar least-squares problems, one per alias
class. For each \(r\), the single weight coefficient \(\widehat\eta_r\) must approximate
the values \(-i\,p\cdot\widehat\epsilon_r\) for \(p\in\Delta_r\), with multiplicities
\(N_p\).

The best common value is just the weighted average of the requests in that class. Define
\begin{equation}
\label{eq:class-mass-center}
N_r
=
\sum_{p\in\Delta_r}N_p,
\qquad
\bar p_r
=
\frac{1}{N_r}
\sum_{p\in\Delta_r}N_p\,p.
\end{equation}
Then the minimizing weight mode is
\begin{equation}
\widehat\eta_r^{\,\star}
=
-i\,\bar p_r\cdot\widehat\epsilon_r.
\label{eq:optimal-class-weight}
\end{equation}
After removing the class average, the residual is
\[
\sum_{p\in\Delta_r}N_p
\bigl|
(p-\bar p_r)\cdot\widehat\epsilon_r
\bigr|^2.
\]
It measures the internal spread of the raw differences within each alias class.

\subsection{The alias-covariance theorem}
The remaining class error has a natural matrix form. Define
\begin{equation}
C_r
=
\sum_{p\in\Delta_r}
N_p\,
(p-\bar p_r)(p-\bar p_r)^T.
\label{eq:alias-covariance-block}
\end{equation}
Then
\[
\sum_{p\in\Delta_r}
N_p\,
\bigl|(p-\bar p_r)\cdot\widehat\epsilon_r\bigr|^2
=
\widehat\epsilon_r^{\,*}C_r\widehat\epsilon_r,
\]
and therefore
\begin{equation}
q(\epsilon)
=
\sum_r
\widehat\epsilon_r^{\,*}C_r\widehat\epsilon_r.
\label{eq:alias-block-decomposition}
\end{equation}
The abstract reduced quadratic form of Section~\ref{sec:node-vs-weight} has become explicit: after
transforming to character modes, it splits into independent positive-semidefinite blocks
\(C_r\), one for each alias class. Each block is simply the weighted covariance, or
scatter matrix, of the raw frequency differences that the sampling grid identifies with
one another.

\begin{theorem}[Alias-covariance block theorem]
\label{thm:alias-blocks}
Let the nodes form a finite subgroup \(\mathcal G\subset\mathbb T^d\) with uniform weights
\(v_k^0=1/S\), and let the rule be exact. Then the finite-group Fourier transform
block-diagonalizes the reduced form on the complexified node tangent space,
\[
\mathcal F\,Q_{\rm red}\,\mathcal F^{-1}=\bigoplus_r C_r ,
\]
the real tangent space being recovered by conjugate pairing, with \(N_r\) and \(\bar p_r\) as in
\eqref{eq:class-mass-center} and \(C_r\) as in \eqref{eq:alias-covariance-block}. For an empty class
\(\Delta_r=\varnothing\) set \(N_r=0\), \(\bar p_r=0\) and \(C_r=0\); the Gram objective
then contains no equation in that class, and \(\widehat\eta_r^\star=0\) selects the
minimum-norm representative. The reduced defect and the optimal weight symbol are
\eqref{eq:alias-block-decomposition} and \eqref{eq:optimal-class-weight}, and the class
minimizers assemble to a real zero-mean weight tangent. Writing
\[
\mathcal A_r=\operatorname{span}\{p-p':p,p'\in\Delta_r\},
\]
one has
\begin{equation}
\label{eq:alias-kernel-rank}
\ker C_r=\mathcal A_r^{\perp},
\qquad
\operatorname{rank}C_r=\dim\mathcal A_r .
\end{equation}
\end{theorem}

\begin{proof}
The classwise minimization \eqref{eq:classwise-frobenius}--\eqref{eq:alias-block-decomposition},
together with Parseval for the reference-weighted inner product \eqref{eq:node-inner-product},
gives the stated Fourier-block decomposition of \(Q_{\rm red}\) and the optimal weight
symbol. Since \(\Delta_{-r}=-\Delta_r\) and \(N_{-p}=N_p\) we have
\(\bar p_{-r}=-\bar p_r\), while reality of \(\epsilon\) gives
\(\widehat\epsilon_{-r}=\overline{\widehat\epsilon_r}\); hence
\(\widehat\eta_{-r}^\star=\overline{\widehat\eta_r^\star}\). In a self-conjugate class
symmetry forces \(\bar p_r=0\), and exactness makes \(\Lambda\) injective modulo
\(L_{\mathcal G}\), so \(\Delta_0=\{0\}\); in both cases \(\widehat\eta_r^\star=0\), and the
class minimizers assemble to a real zero-mean tangent. For \(z\in\mathbb C^d\),
\[
z^*C_rz=\sum_{p\in\Delta_r}N_p\,|(p-\bar p_r)\cdot z|^2 ,
\]
so \(C_r\ge0\) and \(\ker C_r\) is the orthogonal complement of the span of the centered
vectors, which is \(\mathcal A_r\); the rank statement follows.
\end{proof}

\begin{figure}[!b]
\centering
\includegraphics[width=0.82\linewidth]{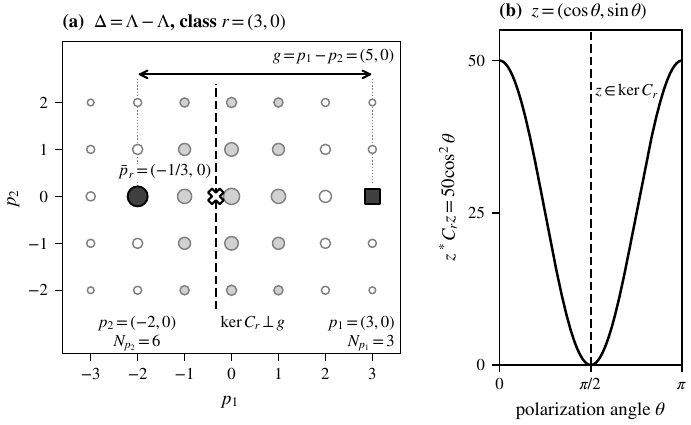}
\caption{Alias-class geometry for the \(4\times3\) band
\(\Lambda=\{0,1,2,3\}\times\{0,1,2\}\) on the \(5\times5\) Cartesian grid. (a) The difference
set \(\Delta=\Lambda-\Lambda\), with marker area proportional to the multiplicity \(N_p\); open
markers are collided differences and filled gray markers are singleton classes. For the residue
\(r=(3,0)\), the class contains \(p_1=(3,0)\) and \(p_2=(-2,0)\) with \(N_{p_1}=3\) and
\(N_{p_2}=6\), center \(\bar p_r=(-1/3,0)\) and alias gap \(g=(5,0)\), so
\(C_r=\frac{ab}{a+b}gg^T=\operatorname{diag}(50,0)\) as in Theorem~\ref{thm:alias-blocks}.
(b) The class defect \(z^*C_rz=50\cos^2\theta\) of a unit polarization
\(z=(\cos\theta,\sin\theta)\). It vanishes only for \(z\perp g\), which spans \(\ker C_r\).}
\label{fig:alias-class-geometry}
\end{figure}

For a class containing just two representatives \(p_1,p_2\), with multiplicities \(a,b\),
the covariance takes the rank-one form
\[
C_r
=
\frac{ab}{a+b}\,gg^T,
\qquad
g=p_1-p_2.
\]
The defect detects only the component of \(\widehat\epsilon_r\) along the alias gap \(g\);
transverse motion remains first-order repairable. Figure~\ref{fig:alias-class-geometry}
illustrates this geometry.

The same decomposition also prices any prescribed linear reweighting rule.
Suppose a compensator acts on character mode \(r\) through a symbol
\(a_r\in\mathbb R^d\),
\begin{equation}
\label{eq:prescribed-symbol}
\widehat\eta_r=-i\,a_r\cdot\widehat\epsilon_r.
\end{equation}
Then its contribution in class \(r\) is
\begin{equation}
\label{eq:prescribed-class-cost}
J_r(a_r;\epsilon)
=
\widehat\epsilon_r^{\,*}C_r\widehat\epsilon_r
+
N_r\left|
(a_r-\bar p_r)\cdot\widehat\epsilon_r
\right|^2 .
\end{equation}
The first term is the irreducible defect; the second is the exact excess paid
for using the prescribed symbol rather than the optimal class center
\(a_r=\bar p_r\). The covariance calculation therefore does more than identify the
optimum: it gives a classwise error formula for any linear-symbol
density-compensation rule.

\begin{proposition}[Prescribed-symbol cost]
\label{prop:prescribed-symbol}
In the setting of Theorem~\ref{thm:alias-blocks}, let real symbols \(a_r\in\mathbb R^d\)
satisfy \(a_{-r}=-a_r\) on non-self-conjugate pairs and \(a_r=0\) on self-conjugate classes.
Then \eqref{eq:prescribed-symbol}, with the complex-bilinear dot product and no conjugation on
\(a_r\), defines a real zero-mean weight tangent whose classwise energy is
\eqref{eq:prescribed-class-cost}.
\end{proposition}

\begin{proof}
The stated conditions give \(\widehat\eta_{-r}=\overline{\widehat\eta_r}\) and
\(\widehat\eta_0=0\), so the tangent is real and zero-mean. In class \(r\) the residual in
\eqref{eq:classwise-frobenius} is \(\sum_{p\in\Delta_r}N_p|(p-a_r)\cdot\widehat\epsilon_r|^2\).
Expanding \(p-a_r=(p-\bar p_r)+(\bar p_r-a_r)\) annihilates the cross term because
\(\sum_{p\in\Delta_r}N_p(p-\bar p_r)=0\), which is \eqref{eq:prescribed-class-cost}.
\end{proof}

For a singleton class the optimal symbol has a familiar geometric
interpretation. If \(\Delta_r=\{p\}\), then \(\bar p_r=p\), so
\begin{equation}
\label{eq:singleton-jacobian-symbol}
\widehat\eta_r^\star
=
-i\,p\cdot\widehat\epsilon_r.
\end{equation}
With the Fourier convention used here, the multiplier
\(-i\,p\cdot\widehat\epsilon_r\) is the discrete divergence symbol associated
with the unique raw representative \(p\) of this singleton class, and not with
the residue label \(r\) that indexes it. Identifying it with the continuous
Fourier coefficient of \(\nabla\!\cdot\epsilon\) requires an appropriate
Fourier interpolant, or more generally an alias-free representation of the
displacement field whose coefficient at the raw frequency \(p\) is
\(\widehat\epsilon_r\). Thus, when the difference set is injective, optimal
first-order reweighting is the sampled infinitesimal Jacobian, or cell-volume,
correction mode by mode \cite{HogeKwanPike1997}. Collisions are precisely where
this picture breaks: several distinct raw differences \(p\) share one character
class, so a single divergence-type weight mode cannot satisfy all of their
requests simultaneously. The optimal response then replaces the individual
\(p\)'s by their multiplicity-weighted center \(\bar p_r\).

\subsection{No collisions}
\label{sec:no-collisions}
When every alias class contains at most one raw difference, every block \(C_r\) vanishes,
and the conclusion extends beyond first order: small node perturbations admit positive
normalized weights with \(G(x,v)=I\) exactly.

\begin{theorem}[Difference-set injective persistence]
\label{thm:injective-persistence}
Let the nodes form a finite subgroup \(\mathcal G\subset\mathbb T^d\) with uniform weights
\(v_k^0=1/S\), and let the rule be exact. The following are equivalent:
\[
q\equiv0,
\qquad
C_r=0\ \text{for every }r,
\qquad
\Delta\longrightarrow\mathbb Z^d/L_{\mathcal G}\ \text{is injective},
\]
equivalently every nonempty \(\Delta_r\) is a singleton. Assume these hold, and let
\(M_R(x)\) be the real matrix whose first row is \(\mathbf 1^T\) and whose remaining rows are
the nonredundant real and imaginary parts of the exactness equations
\(\sum_kv_ke^{ip\cdot x_k}=0\) for nonzero \(p\in\Delta\), with conjugate pairs identified.
With the normalization ordered first, exactness and normalization together read
\[
M_R(x)\,v=b,
\qquad
b=(1,0,\ldots,0)^T .
\]
Then \(M_R(x^0)\) has full row rank, and every node configuration \(x\)
sufficiently near \(x^0\) admits positive normalized weights with \(G(x,v(x))=I\) exactly.
One smooth local choice is
\begin{equation}
\label{eq:persistence-right-inverse}
v(x)=v^0+M_R(x)^T\bigl(M_R(x)M_R(x)^T\bigr)^{-1}\bigl(b-M_R(x)v^0\bigr).
\end{equation}
If \(\Delta\to\mathbb Z^d/L_{\mathcal G}\) is not injective, then \(q\not\equiv0\).
\end{theorem}

\begin{proof}
By Theorem~\ref{thm:alias-blocks}, \(q\equiv0\) if and only if every block \(C_r\) vanishes.
Each \(C_r\) is a scatter matrix with strictly positive multiplicities \(N_p\), so
\(C_r=0\) if and only if all raw vectors in \(\Delta_r\) coincide, that is, every nonempty
\(\Delta_r\) is a singleton; this is exactly injectivity of the reduction map on \(\Delta\).
The last sentence of the theorem is the contrapositive.

Assume injectivity. If some nonzero \(p\in\Delta\) were self-conjugate modulo
\(L_{\mathcal G}\), then \(p\equiv-p\), while symmetry of \(\Delta\) places the two distinct
raw vectors \(p\) and \(-p\) in the same class, contradicting injectivity. Hence no paired
sine row degenerates, and distinct characters give orthogonal rows, so \(M_R(x^0)\) has full
row rank; the normalization row \(\mathbf 1^T\) is independent of them because it is the
character row at \(p=0\), which no nonzero difference meets. Full row rank is an open
condition, so \(M_R(x)M_R(x)^T\) stays invertible nearby and
\eqref{eq:persistence-right-inverse} is well defined and smooth in \(x\); by construction it
satisfies \(M_R(x)v(x)=b\), which is exactness together with normalization. Since
\(v(x^0)=v^0>0\), continuity gives \(v(x)>0\) for \(x\) near \(x^0\).
\end{proof}

Theorem~\ref{thm:injective-persistence} is the Fourier-structured analogue of the openness of strict scalability under a
spanning condition \cite{KutyniokOkoudjouPhilippTuley2013}; compare
Chen--Kutyniok--Okoudjou--Philipp--Wang, Lemma~4.7 \cite{ChenKutyniokOkoudjouPhilippWang2015}.
What is new is that, at a uniform finite-subgroup reference, the relevant rank condition is
characterized exactly by injectivity of \(\Delta\) modulo \(L_{\mathcal G}\), while collisions
are resolved by the covariance blocks \(C_r\). In lattice-rule terms, the hypothesis says that
the sampling lattice is reconstructing for \(\Delta\), equivalently that the reference rule is
exact on \((\Lambda-\Lambda)-(\Lambda-\Lambda)\); Section~\ref{sec:rank1-lattice} shows that
this is strictly stronger than being reconstructing for \(\Lambda\).

At fixed subgroup nodes, changing the exact reference weights does not create
a new obstruction mechanism; it reparametrizes the same one. Let
\(v_k^0=1/S\) be the uniform exact weight and let \(v'\) be any other positive
normalized exact weight on the same nodes. Since the weight-variation space
depends only on the nodes,
\begin{equation}
\label{eq:reference-node-reparametrization}
X_{v'}\epsilon
=
X_{v^0}(S v'\odot\epsilon),
\end{equation}
and therefore
\begin{equation}
\label{eq:reference-q-reparametrization}
q_{v'}(\epsilon)
=
q_{v^0}(S v'\odot\epsilon).
\end{equation}
The value of \(q_{v'}\) on a particular perturbation, and its kernel, do
depend on the chosen exact reference, but the existence of a nontrivial
first-order obstruction does not:
\begin{equation}
\label{eq:reference-obstruction-equivalence}
q_{v'}\equiv0
\quad\Longleftrightarrow\quad
q_{v^0}\equiv0.
\end{equation}
For a fixed perturbation \(\epsilon\), optimizing over the admissible exact
reference weights is consequently a finite convex quadratic program over the
exact-weight polytope.

\section{From harmless directions to exact families}
\label{sec:integrability}

\subsection{The integrability question}
Collisions are different. When some \(C_r\neq 0\), the kernel of \(q\) still identifies
node directions whose first-order defect can be cancelled, but that is only an
infinitesimal statement. A direction \(\epsilon\in\ker q\) need not lie on any actual
nearby family of exact weighted rules. \emph{The next question is geometric
rather than variational: when does a first-order harmless direction integrate to a genuine
local deformation of the exactness set?}

\subsection{Affinely independent collision classes}
The obstruction to integration comes from relations internal to an alias class. For a class
\(\Delta_r\), suppose the representatives are affinely independent: whenever coefficients
\(\lambda_p\) satisfy
\begin{equation}
\label{eq:affine-relation}
\sum_{p\in\Delta_r}\lambda_p=0,
\qquad
\sum_{p\in\Delta_r}\lambda_p\,p=0,
\end{equation}
all \(\lambda_p\) must vanish. The linearized exactness constraints within a class have
coefficient rows \((1,p^T)\), whose independence is exactly this condition, so there are then
no hidden compatibility conditions beyond the first variation and every harmless tangent
direction integrates.

\begin{theorem}[Local exactness manifold under affine independence]
\label{thm:exactness-manifold}
Let \(x^0=\mathcal G\) be finite-subgroup nodes. Assume every nonzero character class is met
by \(\Delta\), equivalently that the normalized exact weight is unique, so that at the
subgroup reference it is the uniform weight \(v^0\). Assume further that every class
\(\Delta_r\) is affinely independent in \(\mathbb R^d\); necessarily
\(|\Delta_r|\le d+1\) in every class, so the hypothesis is restrictive for strongly
colliding higher-dimensional classes. Then the set \(\mathcal E\) of node configurations
admitting exact nonnegative normalized weights is, near \(x^0\), a smooth embedded
submanifold. Nearby exact weights are automatically positive and unique,
\begin{equation}
\label{eq:exactness-tangent-kernel}
T_{x^0}\mathcal E=\ker q,
\end{equation}
and the real codimension of \(\mathcal E\) is
\(\operatorname{rank}Q_{\rm red}\).
\end{theorem}

\begin{proof}
Let \(F(x,v)\) be the nonredundant real exactness map built from
\(\sum_kv_ke^{ip\cdot x_k}=0\) for nonzero \(p\in\Delta\), with conjugate equations paired
and normalization handled by \(\mathbf 1^Tv=1\). By
\eqref{eq:classwise-first-variation}, at the subgroup reference \(DF\) acts classwise as
\((a,z)\mapsto(a+i\,p\cdot z)_{p\in\Delta_r}\), whose row matrix is \((1,ip^T)\); affine
independence of the rows \((1,p^T)\) makes this surjective. To pass to the real paired
equations, take target data with \(h_{-r}=\overline{h_r}\), choose \((a_r,z_r)\) realizing
\(h_r\), and set \(a_{-r}=\overline{a_r}\), \(z_{-r}=\overline{z_r}\); since
\(\Delta_{-r}=-\Delta_r\), the equation at \(-p\) returns \(\overline{a_r+i\,p\cdot z_r}\),
so the conjugate class realizes \(h_{-r}\). If \(r=-r\neq0\) then \(\Delta_r=-\Delta_r\),
and a centrally symmetric affinely independent set contains at most one nonzero pair, so
\(\Delta_r=\{p,-p\}\) and \((\eta,e)\mapsto\eta+i\,p\cdot e\) is onto \(\mathbb C\). The
zero class contributes only normalization because \(\Delta_0=\{0\}\). Hence
\(DF_{(x^0,v^0)}\) is surjective and
\(\mathcal M=\{(x,v):F(x,v)=0,\ \mathbf 1^Tv=1\}\) is a smooth manifold near
\((x^0,v^0)\), with positivity persisting after shrinking the neighborhood.

For \(\pi(x,v)=x\) we have
\(T_{(x^0,v^0)}\mathcal M=\{(\epsilon,\zeta):X\epsilon+Y\zeta=0\}\), and
\(\ker d\pi=\{(0,\zeta):\zeta\in\mathcal Z,\ Y\zeta=0\}=0\) by exact-weight uniqueness, so
\(\pi|_{\mathcal M}\) is an immersion and the constant-rank theorem gives a locally embedded
image. Openness of injectivity of \(Y_x|_{\mathcal Z}\) gives local uniqueness of exact
weights nearby, and compactness of the simplex excludes distant branches: if
\(x_n\to x^0\) carries exact nonnegative normalized \(v_n\), any subsequential limit
\(v^*\) satisfies \(G(x^0,v^*)=I\), hence \(v^*=v^0\). Therefore \(\mathcal E\) is the local
projection of \(\mathcal M\), is embedded, and
\[
T_{x^0}\mathcal E=\{\epsilon:\exists\,\zeta\in\mathcal Z,\ X\epsilon+Y\zeta=0\}=\ker q ,
\]
which is \eqref{eq:exactness-tangent-kernel}. The Fourier description of the kernel and the
codimension count follow from Theorem~\ref{thm:alias-blocks}.
\end{proof}

\begin{remark}[A two-dimensional instance]
\label{rem:exactness-manifold-2d}
The hypotheses of Theorem~\ref{thm:exactness-manifold} are met by the configuration of
Figure~\ref{fig:alias-class-geometry}, \(\Lambda=\{0,\ldots,3\}\times\{0,1,2\}\) on the
\(5\times5\) grid. The difference set \(\Delta=\{-3,\ldots,3\}\times\{-2,\ldots,2\}\) has
\(35\) elements and meets all \(25\) character classes, so the exact weight is unique. Only the
first coordinate collides, at residues \(r_1\in\{2,3\}\), so every class has
\(|\Delta_r|\le2\) and is affinely independent. Hence \(\mathcal E\) is locally a smooth
manifold of real codimension \(\operatorname{rank}Q_{\rm red}=10=|\Delta|-S\) in the
\(50\)-dimensional node space, and \(\dim\ker q=40\); a dense computation of \(Q_{\rm red}\)
confirms both values.
\end{remark}

\subsection{Affine relations}
\label{sec:affine-relations}
When an alias class contains an affine relation, this rank argument breaks. Suppose
coefficients \(\lambda_p\), not all zero, satisfy
\[
\sum_{p\in\Delta_r}\lambda_p=0,
\qquad
\sum_{p\in\Delta_r}\lambda_p\,p=0.
\]
Any constant term and any first-order term linear in \(p\) disappear when combined with these
coefficients. The first variation cannot see this relation. At second order,
however, quadratic terms in \(p\) remain, and they impose an additional compatibility
condition on a supposedly harmless direction. This is the mechanism by which a vector in
\(\ker q\) can fail to integrate to an exact local deformation.

To make that second-order condition explicit, define the matrix-valued Fourier coefficient
\begin{equation}
\label{eq:second-moment-mode}
\widehat M_r(\epsilon)
=
\frac1S\sum_{k=1}^S
\epsilon_k\epsilon_k^T e^{ir\cdot x_k^0}.
\end{equation}
If a first-order harmless direction \(\epsilon\) is tangent to an actual family of exact
configurations, then every affine relation \(\lambda=(\lambda_p)_{p\in\Delta_r}\) in every
class must also satisfy
\begin{equation}
\label{eq:affine-second-order-condition}
\sum_{p\in\Delta_r}
\lambda_p\,p^T\widehat M_r(\epsilon)p
=
0.
\end{equation}
Condition \eqref{eq:affine-second-order-condition} is necessary for exact deformability and
may exclude directions in \(\ker q\) from the tangent cone of the exactness set.

\begin{proposition}[Second-order obstruction from affine relations]
\label{prop:affine-relation-obstruction}
Let \(x^0=\mathcal G\) be uniform finite-subgroup nodes with unique normalized exact weight
\(v^0\). For a class \(r\) set
\[
\mathcal R_r
=
\Bigl\{
\lambda=(\lambda_p)_{p\in\Delta_r}\in\mathbb R^{\Delta_r}:
\sum_{p\in\Delta_r}\lambda_p=0,\
\sum_{p\in\Delta_r}\lambda_p\,p=0
\Bigr\},
\]
and let \(\widehat M_r(\epsilon)\) be as in \eqref{eq:second-moment-mode}. If
\(\epsilon\in T_{x^0}\mathcal E\), then \eqref{eq:affine-second-order-condition} holds for
every class \(r\) and every \(\lambda\in\mathcal R_r\).
\end{proposition}

\begin{proof}
Let \(x_n=x^0+t_n\epsilon+\rho_n\in\mathcal E\) with \(t_n\downarrow0\) and
\(\rho_n=o(t_n)\), and choose exact nonnegative normalized weights \(v_n\). The probability
simplex is compact, so after passing to a subsequence \(v_n\to v^*\), and continuity gives
\(G(x^0,v^*)=I\); exact-weight uniqueness forces \(v^*=v^0\). At \(x^0\) the frame operator
is linear in the weights, so \(G(x^0,v_n)-I=Y(v_n-v^0)\), while Taylor expansion in the node
variable gives
\[
G(x_n,v_n)=G(x^0,v_n)+t_nX_{v_n}\epsilon+o(t_n),
\]
the remainder uniform in \(v_n\) because \(G\) is smooth in \(x\), linear in \(v\), and the
\(v_n\) remain in the compact simplex. Exactness of \(x_n\) therefore gives
\[
0=Y(v_n-v^0)+t_nX_{v_n}\epsilon+o(t_n).
\]
Dividing by \(t_n\) and using \(X_{v_n}\epsilon\to X\epsilon\) shows
\(Y\bigl((v_n-v^0)/t_n\bigr)\to-X\epsilon\). Each \(v_n-v^0\) lies in \(\mathcal Z\) and
\(Y\mathcal Z\) is closed, so \(X\epsilon\in Y\mathcal Z\) and hence \(q(\epsilon)=0\). Hence
there is \(\zeta\in\mathcal Z\) with \(Y\zeta+X\epsilon=0\), unique by injectivity of
\(Y|_{\mathcal Z}\). An injective linear map on the finite-dimensional space \(\mathcal Z\)
has a bounded left inverse; applying it to
\(Y\bigl((v_n-v^0)/t_n\bigr)\to-X\epsilon=Y\zeta\) gives \((v_n-v^0)/t_n\to\zeta\), that is,
\(v_n=v^0+t_n\zeta+\omega_n\) with \(\omega_n=o(t_n)\).

For \(p\in\Delta_r\) write \(M_p(x,v)=\sum_kv_ke^{ip\cdot x_k}\), so \(M_p(x_n,v_n)=0\).
Factoring out \(e^{ip\cdot x^0}\) and expanding to second order, the terms through order
\(t_n^2\) organize by their dependence on \(p\): the direct weight terms, including
\(\omega_n\) at leading order, are constant in \(p\); the terms involving \(\rho_n\) and the
first-order-weight times first-order-node term
\(i t_n^2\,p\cdot\sum_k\zeta_k\epsilon_ke^{ir\cdot x_k^0}\) are linear in \(p\); the
remaining \(\omega_n\) and \(\rho_n\) contributions are \(o(t_n^2)\). The explicit quadratic
node term is \(-\tfrac12t_n^2\,p^T\widehat M_r(\epsilon)p\). Multiplying the moment
equations by \(\lambda_p\) and summing over \(\Delta_r\), the defining identities
\(\sum_p\lambda_p=0\) and \(\sum_p\lambda_pp=0\) annihilate every constant and every
linear term, leaving
\[
0=-\tfrac12t_n^2\sum_{p\in\Delta_r}\lambda_p\,p^T\widehat M_r(\epsilon)p+o(t_n^2).
\]
Dividing by \(t_n^2\) and letting \(n\to\infty\) gives
\eqref{eq:affine-second-order-condition}.
\end{proof}

On product grids, a class colliding in two coordinates contains four representatives forming
a rectangle, and hence a nontrivial affine relation. The two-dimensional tensor DFT provides
a natural test of whether directions in \(\ker q\) integrate to exact rules.

\section{When first-order repairability fails to integrate}
\label{sec:tensor-twist}

\subsection{The first-order kernel}
At the two-dimensional tensor DFT, the first-order kernel is already larger than the set of
obvious exact deformations. For the square tensor DFT the covariance blocks can be written
explicitly. If the grid has shape \(s_1\times s_2\), then for residue \(r=(r_1,r_2)\),
\begin{equation}
\label{eq:tensor-square-covariance}
C_r
=
S\,
\operatorname{diag}\!\bigl(
r_1(s_1-r_1),
\,r_2(s_2-r_2)
\bigr),
\qquad
S=s_1s_2.
\end{equation}
Hence a mode \(\widehat\epsilon_r=(\widehat\epsilon_{r,1},\widehat\epsilon_{r,2})\) lies in
\(\ker C_r\) exactly when each component vanishes whenever the corresponding residue
component is nonzero. Equivalently, the first component may depend only on \(x_2\), while the
second may depend only on \(x_1\). Writing
\begin{equation}
\label{eq:tensor-shear-families}
\mathcal S_1
=
\{(f(x_2),c_2)\},
\qquad
\mathcal S_2
=
\{(c_1,g(x_1))\},
\end{equation}
for the two families this gives
\begin{equation}
\label{eq:tensor-first-order-kernel}
\ker q=\mathcal S_1+\mathcal S_2 .
\end{equation}
Each family is a shear: one may shear in the first coordinate by an arbitrary function of the
second coordinate, or shear in the second coordinate by an arbitrary function of the first,
and each such motion is an exact deformation of the grid. So at this point first-order theory
suggests that a simultaneous row and column shear should also be harmless. The next
calculation shows why that conclusion is false.

\subsection{But the sum does not integrate}
Take a mixed direction \(\epsilon=(f(x_2),g(x_1))\in\mathcal S_1+\mathcal S_2\). A class
colliding in both coordinates has four raw representatives at the corners of a rectangle, and
these satisfy the affine relation \(\lambda=(1,-1,-1,1)\). In the second-order condition
\eqref{eq:affine-second-order-condition} this relation cancels the pure row-shear and pure
column-shear terms, and only the mixed term survives. Two nonconstant shears therefore cannot
be combined into an exact deformation.

\begin{proposition}[Second-order twist obstruction at the tensor DFT]
\label{prop:tensor-twist}
Let \(d=2\), \(\Lambda=\{0,\ldots,s_1-1\}\times\{0,\ldots,s_2-1\}\) and
\(\mathcal G=\mathbb Z_{s_1}\times\mathbb Z_{s_2}\) with \(s_1,s_2\ge2\), uniform weights,
and square tensor nodes \(x^0=\mathcal G\). With \(\mathcal S_1\) and \(\mathcal S_2\) as in
\eqref{eq:tensor-shear-families}, the first-order kernel is
\eqref{eq:tensor-first-order-kernel} and the contingent tangent cone is
\begin{equation}
\label{eq:tensor-tangent-cone}
T_{x^0}\mathcal E
=
\mathcal S_1\cup\mathcal S_2
\subsetneq
\mathcal S_1+\mathcal S_2
=
\ker q .
\end{equation}
Consequently the tangent cone is not a linear space, and the
exactly reweightable node set \(\mathcal E\) is not a smooth embedded manifold at \(x^0\).
\end{proposition}

\begin{proof}
Equation~\eqref{eq:tensor-first-order-kernel} gives \(\ker q=\mathcal S_1+\mathcal S_2\).
Since \(A_0\) is invertible, the exactness equation \(A_0WA_0^*=I\) determines \(W\)
uniquely; Proposition~\ref{prop:exact-square-rigidity} gives \(W=M^{-1}I\). Consequently,
\(T_{x^0}\mathcal E\) is contained in \(\ker q\).

The inclusion \(\mathcal S_1\cup\mathcal S_2\subseteq T_{x^0}\mathcal E\) is exact to all
orders: for \(\epsilon=(f(x_2),c_2)\) take \(x_1(t)=x_1+tf(x_2)\), \(x_2(t)=x_2+tc_2\) with
uniform weights; for a nonzero raw difference \(p\) the \(x_1\)-character sum vanishes when
\(p_1\not\equiv0\pmod{s_1}\), and when \(p_1\equiv0\) exactness reduces to the
\(x_2\)-character sum. The family \(\mathcal S_2\) is symmetric.

Conversely let \(\epsilon=(f(x_2),g(x_1))\in T_{x^0}\mathcal E\). For a residue
\(r=(r_1,r_2)\) with \(r_1,r_2\neq0\) the class contains the four corners
\[
(r_1,r_2),\quad(r_1-s_1,r_2),\quad(r_1,r_2-s_2),\quad(r_1-s_1,r_2-s_2),
\]
on which \(\lambda=(1,-1,-1,1)\) is an affine relation. Writing
\(\operatorname{Tw}_r[h]\) for the corresponding alternating sum, one has
\(\operatorname{Tw}_r[a+b\cdot p]=0\) for every affine function of \(p\), while
\(\operatorname{Tw}_r[p_1p_2]=s_1s_2\). For \(\epsilon=(f(x_2),g(x_1))\) the pure \(f^2\)
term vanishes after summing over \(x_1\) and the pure \(g^2\) term after summing over
\(x_2\), so the mixed part of \(p^T\widehat M_r(\epsilon)p\) is
\(2p_1p_2\,\widehat g_{r_1}\widehat f_{r_2}\). Proposition~\ref{prop:affine-relation-obstruction}
therefore gives
\begin{equation}
\label{eq:tensor-product-condition}
0=2\,\operatorname{Tw}_r[p_1p_2]\,\widehat g_{r_1}\widehat f_{r_2}
=2s_1s_2\,\widehat g_{r_1}\widehat f_{r_2}.
\end{equation}
Since this holds for every nonzero residue
pair, either all nonzero Fourier coefficients of \(f\) vanish or all those of \(g\) do, so
\(\epsilon\in\mathcal S_1\cup\mathcal S_2\). Taking \(f\) and \(g\) both nonconstant gives an
element of \(\mathcal S_1+\mathcal S_2\) outside the union, so the inclusion is strict. A
union of two crossing linear branches is not a linear space, hence cannot be the tangent
space of a smooth embedded manifold.
\end{proof}

\subsection{Size of the obstruction}
For \(\epsilon\in\ker q\), the first-order Gram defect can be eliminated, but a
second-order residual may remain. We quantify its smallest attainable squared coefficient by
\begin{equation}
\label{eq:intrinsic-second-order-defect}
\Phi_*(\epsilon)
=
\inf
\liminf_{t\downarrow0}
\frac{\|G(x(t),v(t))-I\|_F^2}{t^4},
\end{equation}
where the infimum is over admissible paths
\[
x(t)=x^0+t\epsilon+o(t),
\qquad
v(t)\ge0,
\qquad
\mathbf 1^Tv(t)=1,
\qquad
v(t)\to v^0 .
\]

\begin{figure}[tbp]
\centering
\includegraphics[width=0.9\linewidth]{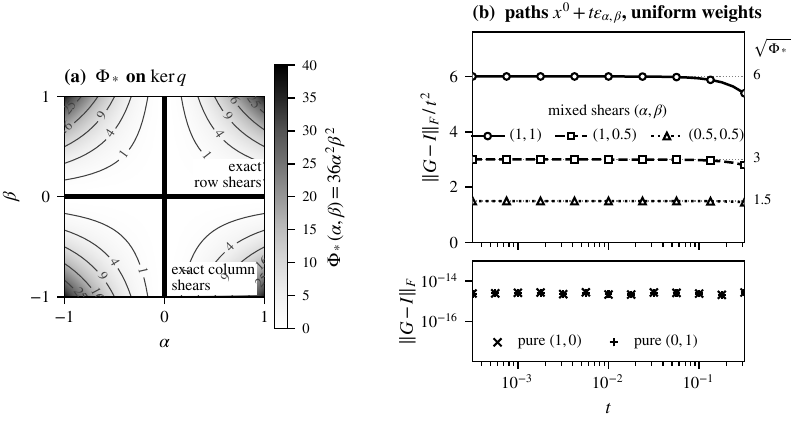}
\caption{First-order kernel versus nonlinear integrability at the \(4\times4\) tensor DFT, in
the family \(\epsilon_{\alpha,\beta}=(\alpha\cos x_2,\beta\cos x_1)\subset\ker q\). (a) The
optimal second-order defect \(\Phi_*(\alpha,\beta)=36\alpha^2\beta^2\) of
Proposition~\ref{prop:tensor-second-order-factorization}. The coordinate axes, drawn bold, are
the exact row and column shears of Proposition~\ref{prop:tensor-twist}. (b) Straight paths
\(x^0+t\epsilon_{\alpha,\beta}\) with uniform weights. For mixed shears,
\(\|G-I\|_F/t^2\) tends to \(\sqrt{\Phi_*}=6|\alpha\beta|\) (dotted levels). Because pure shears
are exact, the \(t^2\) term of a mixed path is the cross term alone, and
Proposition~\ref{prop:tensor-second-order-factorization} shows that no admissible correction
reduces it. For the pure shears \((1,0)\) and \((0,1)\), the lower panel shows that
\(\|G-I\|_F\) itself stays at roundoff for every \(t\).}
\label{fig:tensor-integrability}
\end{figure}

For a mixed tensor shear \(\epsilon=(f(x_2),g(x_1))\), the second-order problem also
separates cleanly. Define the one-dimensional square-DFT energy
\begin{equation}
\label{eq:square-dft-energy}
E_s(h)
=
s\sum_{r=1}^{s-1} r(s-r)\,|\widehat h_r|^2 .
\end{equation}
After optimizing over all admissible second-order corrections, the residual is the product of
the two one-dimensional energies (Proposition~\ref{prop:tensor-second-order-factorization}), so
the obstruction is genuinely mixed and vanishes whenever either shear is constant.
Figure~\ref{fig:tensor-integrability} shows it over a two-parameter cosine family, together with
the second-order behavior of straight paths.

\begin{proposition}[Tensor factorization of the optimal second-order obstruction]
\label{prop:tensor-second-order-factorization}
For the square tensor DFT and every real \(q\)-null tangent
\(\epsilon=(f(x_2),g(x_1))\), let \(\Phi_*(\epsilon)\) be as in
\eqref{eq:intrinsic-second-order-defect}, the infimum being over the admissible path class
defined there. Let \(R_2(\epsilon;\xi,\omega)\) denote the coefficient of \(t^2\) in
\(G(x^0+t\epsilon+t^2\xi,\ v^0+t^2\omega)-I\), and set
\[
\Phi_{\rm opt}(f,g)
=
\min_{\xi,\ \omega:\ \mathbf 1^T\omega=0}
\|R_2(\epsilon;\xi,\omega)\|_F^2 .
\]
Then \(\Phi_*(\epsilon)=\Phi_{\rm opt}(f,g)\), and with \(E_s\) as in
\eqref{eq:square-dft-energy},
\begin{equation}
\label{eq:tensor-second-order-factorization}
\Phi_*(\epsilon)
=
E_{s_1}(g)\,E_{s_2}(f).
\end{equation}
In particular, for \(s_1=s_2=4\) and
\(\epsilon_{\alpha,\beta}=(\alpha\cos x_2,\beta\cos x_1)\),
\begin{equation}
\label{eq:tensor-4x4-cosine}
\Phi_*(\alpha,\beta)
=
36\alpha^2\beta^2 .
\end{equation}
\end{proposition}

\begin{proof}
At the square reference \(X\epsilon=0\) for these \(q\)-null tangents and the exact
normalized weight is unique. Write \(x(t)=x^0+t\epsilon+r(t)\) with \(r(t)=o(t)\) and
\(v(t)=v^0+u(t)\) with \(\mathbf 1^Tu(t)=0\) and \(u(t)\to0\). Along any sequence on which the
second-order quotient is bounded, \(G-I=O(t^2)\); since \(G\) is linear in \(v\) and
\(X\epsilon=0\), the first-order terms give \(Yu=-Xr+O(t^2)+o(|u|)\). Injectivity of
\(Y|_{\mathcal Z}\) then gives \(u=O(|r|+t^2)=o(t)\). Set \(z(t)=(r(t),u(t))\) and
\(L(r,u)=Xr+Yu\). Taylor expansion gives
\[
G(x(t),v(t))-I=Lz(t)+t^2B(\epsilon,\epsilon)+o(t^2),
\]
since the \(t\epsilon\)-by-\(z(t)\) terms and the terms quadratic in \(z(t)\) are
\(o(t^2)\). Decompose \(z=z_0+z_\perp\) with \(z_0\in\ker L\) and
\(z_\perp\perp\ker L\). As \(L\) is bounded below on \((\ker L)^\perp\), a finite
second-order quotient forces \(z_\perp=O(t^2)\), while the kernel part has zero linear image
and interacts with \(t\epsilon\) only at order \(o(t^2)\). Along a subsequence realizing the
liminf, a convergent subsequence of \(z_\perp/t^2\) supplies an effective correction
\((\xi,\omega)\), so the limiting coefficient is at least \(\Phi_{\rm opt}\). Conversely each
fixed \((\xi,\omega)\) with \(\mathbf 1^T\omega=0\) is realized by
\((x^0+t\epsilon+t^2\xi,\ v^0+t^2\omega)\), which is positive for small \(t\) because
\(v^0_k=1/S>0\). Hence \(\Phi_*=\Phi_{\rm opt}\).

Only classes with \(r_1,r_2\neq0\) contribute. Their representatives are
\(p_1\in\{r_1,r_1-s_1\}\) and \(p_2\in\{r_2,r_2-s_2\}\) with product multiplicities
\(N_p=N^{(1)}_{p_1}N^{(2)}_{p_2}\), where \(N^{(a)}_{r_a}=s_a-r_a\) and
\(N^{(a)}_{r_a-s_a}=r_a\); the weighted coordinate means vanish. The corrections \(\xi\) and
\(\omega\) contribute, on each class, functions in \(\operatorname{span}\{1,p_1,p_2\}\), and
the mixed monomial is orthogonal to that span because
\[
\sum_pN_pp_1p_2=\sum_pN_pp_1^2p_2=\sum_pN_pp_1p_2^2=0,
\]
each identity factoring through a vanishing first moment in one coordinate. The mixed
quadratic term is therefore already the least-squares residual after optimizing all
second-order corrections. For \(\epsilon=(f(x_2),g(x_1))\) the pure \(f^2\) and \(g^2\)
terms vanish in these classes by character orthogonality, and the mixed part of
\(-\tfrac12(p\cdot\epsilon)^2\) is \(-p_1p_2\,\widehat g_{r_1}\widehat f_{r_2}\). The
classwise optimal squared residual is thus
\[
|\widehat g_{r_1}|^2|\widehat f_{r_2}|^2
\Bigl(\sum_{p_1}N^{(1)}_{p_1}p_1^2\Bigr)
\Bigl(\sum_{p_2}N^{(2)}_{p_2}p_2^2\Bigr),
\]
and the one-dimensional square identity
\(\sum_{p_a\in\Delta^{(a)}_{r_a}}N^{(a)}_{p_a}p_a^2=s_ar_a(s_a-r_a)\) makes the sum over
\(r_1,r_2\neq0\) factor as \(E_{s_1}(g)E_{s_2}(f)\). For \(s=4\) one has
\(\widehat{\cos}_1=\widehat{\cos}_3=1/2\), so
\(E_4(\alpha\cos)=4(1\cdot3+3\cdot1)\alpha^2/4=6\alpha^2\) and likewise
\(E_4(\beta\cos)=6\beta^2\), giving \eqref{eq:tensor-4x4-cosine}.
\end{proof}

\section{One dimension: what oversampling buys}
\label{sec:one-dim}

\subsection{The collision window}
For one dimension, the class geometry can be read off directly. Take \(m\) consecutive
Fourier modes and \(s\) equally spaced sample points, with
\begin{equation}
\label{eq:1d-deficit}
m\le s\le 2m-2,
\qquad
\delta:=2m-s.
\end{equation}
The raw differences are
\[
\Delta=\{-(m-1),\ldots,m-1\},
\qquad
N_p=m-|p|.
\]
Only the differences near the two ends of this interval collide modulo \(s\). As a concrete
example, let \(m=8\) and \(s=10\), so \(\delta=6\). The five collided classes are
\[
\begin{array}{c|cc|cc}
r & p^+ & p^- & N_{p^+} & N_{p^-}\\
\hline
3 & 3 & -7 & 5 & 1\\
4 & 4 & -6 & 4 & 2\\
5 & 5 & -5 & 3 & 3\\
6 & 6 & -4 & 2 & 4\\
7 & 7 & -3 & 1 & 5
\end{array}
\]
while every other class contains a single raw difference and is exactly compensable
to first order.

\begin{figure}[tbp]
\centering
\includegraphics[width=0.9\linewidth]{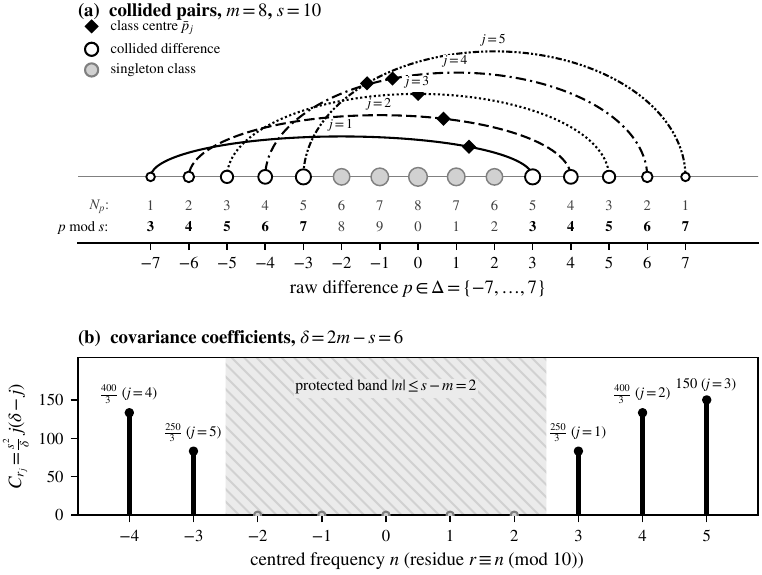}
\caption{The one-dimensional collision window for \(m=8\), \(s=10\), \(\delta=2m-s=6\). (a) The
raw differences \(p\in\{-7,\ldots,7\}\), with marker area proportional to \(N_p=m-|p|\) and
residue \(p\bmod s\) listed below. The arcs join the five collided pairs
\((p_j^+,p_j^-)=(3,-7),\ldots,(7,-3)\), and the diamonds mark the class centers
\(\bar p_j=4/3,2/3,0,-2/3,-4/3\). (b) The covariance coefficients
\(C_{r_j}=(s^2/\delta)\,j(\delta-j)\) of Theorem~\ref{thm:deficit} plotted against centered
frequency; residue \(5\) is the self-conjugate Nyquist class. Every class in the protected band
\(|n|\le s-m=2\) is a singleton with \(C_r=0\).}
\label{fig:collision-window}
\end{figure}

The example makes the general pattern transparent. The collided residues are
\begin{equation}
\label{eq:1d-collided-residues}
r_j=s-m+j,
\qquad
j=1,\ldots,\delta-1,
\end{equation}
and each such class contains exactly two representatives,
\begin{equation}
\label{eq:1d-collision-representatives}
p_j^+=r_j,
\qquad
p_j^-=r_j-s=j-m.
\end{equation}
Their multiplicities are
\begin{equation}
\label{eq:1d-collision-multiplicities}
N_{p_j^+}=\delta-j,
\qquad
N_{p_j^-}=j.
\end{equation}
Thus there are exactly
\[
\delta-1=2m-s-1
\]
collided classes. As \(s\) increases toward the difference-set injective threshold \(2m-1\),
this collision window shrinks one class at a time; at the threshold it disappears altogether. Figure~\ref{fig:collision-window} shows the window for
\(m=8\), \(s=10\).

\subsection{The deficit law}
For these two-point classes, the optimal response follows from the general class-average
formula. Using the multiplicities above,
\begin{equation}
\label{eq:1d-alias-center}
\bar p_j
=
\frac{(\delta-j)p_j^+ + j p_j^-}{\delta}
=
\frac{(m-\delta)(\delta-2j)}{\delta},
\end{equation}
so
\begin{equation}
\label{eq:1d-optimal-weight-response}
\widehat\eta_{r_j}^{\,\star}
=
-i\,\bar p_j\,\widehat\epsilon_{r_j}.
\end{equation}
The center \(\bar p_j\) records the balance between the two aliased representatives. Near the
middle of the collision window their multiplicities are nearly equal and the optimal weight
correction is small; when \(\delta\) is even, the central class has \(\bar p_j=0\). At the
square endpoint \(s=m\), so \(\delta=m\), every class center vanishes and the optimal
first-order weight correction is zero throughout.

For a two-point class, the residual is just the weighted variance of its two representatives.
Using
\begin{equation}
\label{eq:1d-alias-gap}
p_j^+-p_j^-=s
\end{equation}
and the multiplicities \(\delta-j\) and \(j\), we obtain
\begin{equation}
\label{eq:1d-class-residual}
\sum_{p\in\Delta_{r_j}}
N_p\,
\bigl|(p-\bar p_j)\widehat\epsilon_{r_j}\bigr|^2
=
\frac{s^2}{\delta}\,
j(\delta-j)\,
|\widehat\epsilon_{r_j}|^2.
\end{equation}
Summing over the \(\delta-1\) collided classes gives
\begin{equation}
\label{eq:1d-deficit-law}
q(\epsilon)
=
\frac{s^2}{\delta}
\sum_{j=1}^{\delta-1}
j(\delta-j)\,
|\widehat\epsilon_{r_j}|^2.
\end{equation}
Every non-collided class contributes zero. The one-dimensional obstruction is
concentrated entirely in the collision window, with a parabolic weight \(j(\delta-j)\) across
that window.

\begin{theorem}[One-dimensional deficit law]
\label{thm:deficit}
For consecutive \(m\)-mode sampling with \(m\le s\le2m-2\), every class outside the collided
residues \eqref{eq:1d-collided-residues} is a singleton. The collided class \(r_j\) has center
\(\bar p_{r_j}\) given by \eqref{eq:1d-alias-center} and covariance
\begin{equation}
\label{eq:1d-class-covariance}
C_{r_j}=\frac{s^2}{\delta}\,j(\delta-j),
\end{equation}
and the reduced defect is \eqref{eq:1d-deficit-law}.
\end{theorem}

\begin{proof}
The representatives and multiplicities \eqref{eq:1d-collision-representatives}--%
\eqref{eq:1d-collision-multiplicities} give the center \eqref{eq:1d-alias-center}, and the
two-point variance calculation \eqref{eq:1d-class-residual} gives
\eqref{eq:1d-class-covariance}. All remaining classes are singletons and contribute zero;
summing by Theorem~\ref{thm:alias-blocks} gives \eqref{eq:1d-deficit-law}.
\end{proof}

\subsection{Which motions are repairable?}
Equation~\eqref{eq:1d-deficit-law} makes the repairable directions immediate. Since every
coefficient in the sum is positive,
\[
q(\epsilon)=0
\]
if and only if the Fourier coefficients in every collided class vanish. Using centered
discrete frequencies, the collided residues are precisely those outside the band
\(|n|\le s-m\), so
\begin{equation}
\label{eq:1d-repairable-band}
q(\epsilon)=0
\qquad\Longleftrightarrow\qquad
\operatorname{supp}\widehat\epsilon
\subseteq
\{|n|\le s-m\}.
\end{equation}
Oversampling creates a protected low-frequency band of node motions that can be
cancelled completely by reweighting. At the square endpoint \(s=m\), only the rigid mode
\(n=0\) is protected; as \(s\) increases, the protected band widens until the collision window
disappears at \(s=2m-1\).

This low-frequency tangent space has a classical counterpart in the theory of orthogonal
polynomials on the unit circle (OPUC), specifically in Szeg\H{o} quadrature
\cite{JonesNjastadThron1989,Simon2005OPUC1}. After shifting the band to
\(\{0,\ldots,m-1\}\), an exact weighted rule is a positive probability measure
\(\mu=\sum_{k=1}^sv_k\delta_{e^{ix_k}}\) whose moments \(\widehat\mu(n)\) vanish for
\(n=1,\ldots,m-1\) \cite{Peherstorfer2011}. Positive \(s\)-atomic probability measures on the
circle, with node sets taken up to ordering, are parametrized by their Verblunsky coefficients
\[
\alpha_0,\ldots,\alpha_{s-2}\in\mathbb D,
\qquad
\alpha_{s-1}\in\partial\mathbb D,
\]
and the moment conditions are equivalent to
\(\alpha_0=\cdots=\alpha_{m-2}=0\). The free parameters
\(\alpha_{m-1},\ldots,\alpha_{s-2}\in\mathbb D\) and \(\alpha_{s-1}\in\partial\mathbb D\) give
the exact-rule family real dimension
\[
2(s-m)+1,
\]
which agrees with the dimension of the band \(|n|\le s-m\). The existence and dimension of
the one-dimensional exactness manifold are therefore classical; what the present calculation
adds is an explicit Fourier-coordinate description of its tangent directions together with the
transverse defect \(q\), which quantifies what happens when a perturbation leaves that
manifold.

\begin{corollary}[Characterization of first-order repairable modes]
\label{cor:low-frequency-compensable}
In the one-dimensional deficit family \(m\le s\le2m-2\), a real node tangent is first-order
repairable if and only if its centered discrete Fourier support satisfies
\eqref{eq:1d-repairable-band}. Moreover \(T_{x^0}\mathcal E=\ker q\), and this space has real
dimension \(2(s-m)+1\).
\end{corollary}

\begin{proof}
By \eqref{eq:1d-collided-residues}, the collided classes occupy precisely the centered
frequencies \(|n|>s-m\). Their covariance coefficients satisfy
\[
C_{r_j}=\frac{s^2}{\delta}j(\delta-j)>0,
\]
so Theorem~\ref{thm:deficit} gives \eqref{eq:1d-repairable-band}. Every collision class
contains at most two distinct representatives and is therefore affinely independent; the
consecutive difference set \(\Delta\), of cardinality \(2m-1>s\), meets every nonzero
character class. Theorem~\ref{thm:exactness-manifold} then gives
\(T_{x^0}\mathcal E=\ker q\). The surviving Fourier modes comprise one real zero mode and
\(s-m\) complex conjugate pairs, yielding real dimension \(2(s-m)+1\).
\end{proof}

\subsection{Random jitter}
The covariance blocks also describe the mean-square response to correlated node motion. Let
\(\epsilon:\mathcal G\to\mathbb R^d\) be a mean-zero, second-order stationary random field on
the sampling group, with covariance \(K(h)=\mathbb E[\epsilon(x+h)\epsilon(x)^T]\) and spectral
covariance matrices
\begin{equation}
\label{eq:stationary-spectral-covariance}
\mathcal S_r
:=\mathbb E\bigl[\widehat\epsilon_r\widehat\epsilon_r^{\,*}\bigr]
=\frac1S\sum_{h\in\mathcal G}K(h)e^{ir\cdot h}.
\end{equation}
The fixed-weight defect involves the uncentered second moments
\begin{equation}
\label{eq:uncentered-class-moment}
D_r
:=\sum_{p\in\Delta_r}N_p\,pp^T
=C_r+N_r\bar p_r\bar p_r^T .
\end{equation}

\begin{proposition}[Stationary jitter]
\label{prop:stationary-jitter}
In the setting of Theorem~\ref{thm:alias-blocks}, let \(\epsilon\) be a mean-zero
second-order stationary field with spectral covariances
\eqref{eq:stationary-spectral-covariance}. Then
\[
\mathbb E\,q(\epsilon)=\sum_r\operatorname{tr}(C_r\mathcal S_r),
\qquad
\mathbb E\|X\epsilon\|_F^2=\sum_r\operatorname{tr}(D_r\mathcal S_r),
\]
and \(\mathbb E\,q(\epsilon)=0\) if and only if
\(\operatorname{ran}\mathcal S_r\subseteq\ker C_r\) for every \(r\). Whenever the second
expectation is positive, the fraction of expected fixed-weight energy remaining after optimal
reweighting,
\[
\rho_{\mathrm{stat}}
=
\frac{\sum_r\operatorname{tr}(C_r\mathcal S_r)}{\sum_r\operatorname{tr}(D_r\mathcal S_r)},
\]
satisfies \(0\le\rho_{\mathrm{stat}}\le1\).
\end{proposition}

\begin{proof}
Taking expectations in \eqref{eq:alias-block-decomposition}, and in
\eqref{eq:classwise-frobenius} with \(\widehat\eta=0\), gives the two identities, since
\(\mathbb E[\widehat\epsilon_r^{\,*}B\widehat\epsilon_r]=\operatorname{tr}(B\mathcal S_r)\)
for any matrix \(B\). Both \(C_r\) and \(\mathcal S_r\) are positive semidefinite, so
\(\operatorname{tr}(C_r\mathcal S_r)=\|C_r^{1/2}\mathcal S_r^{1/2}\|_F^2\ge0\), with
equality exactly when \(\operatorname{ran}\mathcal S_r\subseteq\ker C_r\). Finally \(0\le C_r\le D_r\) gives
\(0\le\rho_{\mathrm{stat}}\le1\).
\end{proof}

In the one-dimensional deficit family the protected band \(|n|\le s-m\) has \(C_r=0\), so
jitter supported spectrally in this band has no irreducible first-order defect. Compensation
depends on how the perturbation spectrum overlaps the damaging alias classes, rather than on
correlation alone. For independent, mean-zero node jitters with common covariance \(\Gamma\),
stationarity gives \(K(h)=\Gamma\mathbf1_{\{h=0\}}\), and the preceding identities reduce to
\begin{equation}
\label{eq:jitter-mode-covariance}
\mathbb E\!\left[
  \widehat\epsilon_r\widehat\epsilon_r^{\,*}
\right]
=\frac{\Gamma}{S},
\qquad
\mathbb E\,q(\epsilon)
=\frac1S\sum_r\operatorname{tr}(C_r\Gamma).
\end{equation}
In one dimension this gives a closed form for the repairable fraction.

\begin{corollary}[Mean-square jitter repair fraction]
\label{cor:jitter-repair}
Let the node jitter be iid real with mean zero and variance \(\sigma^2\). In the
one-dimensional deficit family with \(m\le s\le2m-1\),
\begin{equation}
\label{eq:1d-jitter-expected-defect}
\mathbb E\,q(\epsilon)=\frac{\sigma^2s}{6}\bigl(\delta^2-1\bigr),
\qquad
\mathbb E\|H^{(0)}\|_F^2
=\frac{\sigma^2}{s}\sum_{p\ne0}(m-|p|)p^2
=\frac{\sigma^2}{6s}m^2(m^2-1),
\end{equation}
where \(H^{(0)}\) is the first-order Gram perturbation at fixed weights. Hence the fraction
of expected fixed-weight first-order Gram damage remaining after optimal reweighting is
\begin{equation}
\label{eq:1d-jitter-rho}
\rho_{m,s}
=
\frac{s^2(\delta^2-1)}
     {m^2(m^2-1)},
\end{equation}
the repairable fraction is \(1-\rho_{m,s}\), and in the proportional limit
\(s/m\to\kappa\in[1,2]\),
\begin{equation}
\label{eq:1d-jitter-limit}
\rho_{m,s}
\longrightarrow
\kappa^2(2-\kappa)^2,
\qquad
1-\rho_{m,s}
\longrightarrow
1-\kappa^2(2-\kappa)^2 .
\end{equation}
These are ratios of expected quadratic energies, not expectations of realization-wise ratios.
\end{corollary}

\begin{proof}
Since \(\mathbb E|\widehat\epsilon_r|^2=\sigma^2/s\) for every nonzero residue,
\eqref{eq:1d-deficit-law} gives
\[
\mathbb E\,q(\epsilon)
=\frac{s^2}{\delta}\sum_{j=1}^{\delta-1}j(\delta-j)\frac{\sigma^2}{s}
=\frac{\sigma^2s}{6}\bigl(\delta^2-1\bigr).
\]
At fixed weights \(\|H^{(0)}\|_F^2=\sum_{p\ne0}(m-|p|)p^2|\widehat\epsilon_p|^2\), and
\(\sum_{p\ne0}(m-|p|)p^2=m^2(m^2-1)/6\) gives the second identity. The ratio is
\eqref{eq:1d-jitter-rho}, and substituting \(\delta=2m-s\) with \(s/m\to\kappa\) gives
\(\delta/m\to2-\kappa\) and hence \eqref{eq:1d-jitter-limit}.
\end{proof}

\begin{figure}[tbp]
\centering
\includegraphics[width=0.66\linewidth]{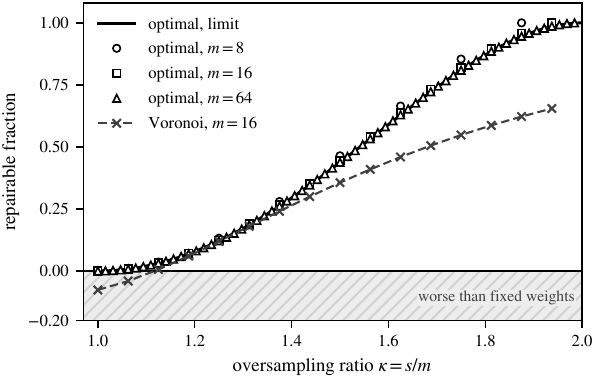}
\caption{Mean-square first-order payoff in the one-dimensional deficit family. The vertical
axis is the fraction of expected fixed-weight first-order Gram energy removed by reweighting;
zero means no improvement, and the hatched region means the reweighting makes the error worse.
The open markers are the exact finite values of Corollary~\ref{cor:jitter-repair} for
\(m=8,16,64\), and the solid curve is its proportional limit. The dashed curve is the same observable
for the prescribed Voronoi weights at \(m=16\); it is negative at critical sampling.}
\label{fig:repair}
\end{figure}

Figure~\ref{fig:repair} compares the exact finite values with the limit and with the
prescribed Voronoi rule of the next subsection.

\subsection{Why Voronoi misses the optimum}
Voronoi weights are a classical geometric density-compensation rule in irregular sampling
\cite{FeichtingerGroechenigStrohmer1995}. Here they provide a useful comparison because their
first variation is fixed by the local cell geometry rather than chosen to minimize the Gram
defect. For ordered periodic nodes,
\begin{equation}
\label{eq:voronoi-weights}
v_k(t)
=
\frac{x_{k+1}(t)-x_{k-1}(t)}{4\pi},
\qquad
\eta_k
=
\frac{s}{4\pi}
\bigl(\epsilon_{k+1}-\epsilon_{k-1}\bigr),
\end{equation}
and the discrete Fourier transform gives
\begin{equation}
\label{eq:voronoi-symbol}
\widehat\eta_r
=
-i\,a_r^{\rm Vor}\widehat\epsilon_r,
\qquad
a_r^{\rm Vor}
=
\frac{s}{2\pi}
\sin\!\left(\frac{2\pi r}{s}\right).
\end{equation}
Voronoi therefore replaces the wavenumber demanded by the class geometry with the
centered-difference modified wavenumber \(a_r^{\rm Vor}\). In an injective class the optimum
is the raw representative \(p\); in a collided class it is the alias center \(\bar p_r\) of
\eqref{eq:1d-alias-center}. By Proposition~\ref{prop:prescribed-symbol}, the excess over the
Frobenius optimum per unit modal amplitude \(|\widehat\epsilon_r|^2\) is
\begin{equation}
\label{eq:voronoi-excess}
E_r^{\rm Vor}
=
N_r\,
\bigl|a_r^{\rm Vor}-\bar p_r\bigr|^2 ,
\end{equation}
so Voronoi is optimal in a class precisely when its symbol lands on the class center. For the
collided one-dimensional classes \(r_j=s-m+j\), where \(N_{r_j}=\delta\),
\begin{equation}
\label{eq:voronoi-deficit-excess}
E_{r_j}^{\rm Vor}
=
\delta
\left[
\frac{s}{2\pi}
\sin\!\left(\frac{2\pi r_j}{s}\right)
-
\frac{(m-\delta)(\delta-2j)}{\delta}
\right]^2 .
\end{equation}
Voronoi can nevertheless be accurate on low modes, because
\(a_r^{\rm Vor}=r+O(r^3/s^2)\) for \(|r|\ll s\); geometric cell-volume weighting tracks the
optimum for smooth, low-frequency node motion even though it is not classwise optimal in
general. The distinction is clearest at \(s=2m-1\): the difference set is injective, so optimal
reweighting cancels every first-order node perturbation exactly, whereas
\(a_r^{\rm Vor}=r\) only for the rigid mode \(r=0\). Complete first-order repairability by
optimized weights does not imply exactness of the prescribed Voronoi rule.

At critical sampling \(s=m\), the optimal first-order weight correction vanishes, whereas
Voronoi generally introduces a nonzero one (it vanishes only when the perturbation is supported
on the mode \(r=0\) and, for even \(s\), \(r=s/2\), where \(a^{\rm Vor}_r=0\)). By the square-rule orthogonality established in
Section~\ref{sec:critical}, this increases the Frobenius defect relative to fixed weights.
The negative Voronoi values in Figure~\ref{fig:repair} illustrate the consequence: prescribed
geometric weighting can be worse than doing nothing.

\section{Critical sampling: why weights lose leverage}
\label{sec:critical}

\subsection{First-order reweighting orthogonality}
At critical sampling the loss of reweighting freedom is stronger than the one-dimensional
calculation suggests: weights lose all first-order leverage, in every dimension.

\begin{proposition}[Square-rule reweighting orthogonality]
\label{prop:exact-square-rigidity}
Let \(M=S\) and suppose
\begin{equation}
\label{eq:square-exact-rule}
A_0W_0A_0^*=I
\end{equation}
with \(W_0=\operatorname{diag}(v^0)\), \(v_k^0>0\) and \(\sum_kv_k^0=1\). Then \(A_0\) is
invertible, the weights are uniform,
\begin{equation}
\label{eq:square-uniform-weight}
v_k^0=\frac1M,
\end{equation}
\(A_0/\sqrt M\) is unitary (equivalently \(A_0\) is a complex Hadamard matrix), and the
node-generated and weight-generated first-order Gram perturbations are Frobenius-orthogonal,
\begin{equation}
\label{eq:square-tangent-orthogonality}
\langle X\epsilon,Y\zeta\rangle_F=0
\qquad
\text{for every }\epsilon,\zeta,
\end{equation}
or equivalently
\begin{equation}
\label{eq:square-reweighting-orthogonality}
P_YX=0.
\end{equation}
Thus diagonal reweighting cannot reduce the first-order Frobenius node defect at any positive
exact square Fourier rule, and the best first-order correction is to leave the weights
unchanged.
\end{proposition}

\begin{proof}
Since \(A_0W_0A_0^*=I\) is square with \(W_0>0\), the matrix \(B_0=A_0W_0^{1/2}\) is
unitary, so
\[
1=(B_0^*B_0)_{kk}=\sum_{\ell\in\Lambda}|(A_0)_{\ell k}|^2v_k^0=Mv_k^0,
\]
because every Fourier entry has modulus one; this is
\eqref{eq:square-uniform-weight}, and \(A_0/\sqrt M\) is then unitary. Write
\(\Omega_a=\operatorname{diag}(\ell_a)_{\ell\in\Lambda}\). The \(a\)-coordinate node
derivative is \(i[\Omega_a,B_a]\) with
\(B_a=A_0\operatorname{diag}(v^0\epsilon_a)A_0^*\), while a weight tangent is
\(Y(\zeta)=A_0\operatorname{diag}(\zeta)A_0^*\). Both \(B_a\) and \(Y(\zeta)\) lie in the
commuting algebra \(A_0(\text{diagonal})A_0^{-1}\), so
\[
\langle i[\Omega_a,B_a],Y(\zeta)\rangle_F
=\operatorname{Re}\,i\operatorname{tr}\bigl(\Omega_a[B_a,Y(\zeta)]\bigr)=0 .
\]
Summing over \(a\) gives \eqref{eq:square-tangent-orthogonality}, and hence
\eqref{eq:square-reweighting-orthogonality}.
\end{proof}

\begin{remark}[Determinant interpretation]
\label{rem:square-determinant}
Along every normalized path \(\operatorname{tr}G=\sum_kv_k\sum_\ell|A_{\ell k}|^2=M\), so
\(G=I\) maximizes \(\det G\) at fixed trace by the arithmetic--geometric mean inequality. This
is the determinant criterion of \(D\)-optimal design and the setting of the Kiefer--Wolfowitz
equivalence theory \cite{KieferWolfowitz1960}; here it gives a second route to
\eqref{eq:square-tangent-orthogonality}. Since \(A\) is square,
\begin{equation}
\label{eq:square-logdet-separation}
\log\det G
=
\log|\det A(x)|^2+\sum_k\log v_k .
\end{equation}
(i) The right-hand side is a sum of a function of \(x\) and a function of \(v\), so its Hessian
at \((x^0,v^0)\) has no mixed node-weight block. (ii) Along a normalized straight path
\((x^0+t\epsilon,v^0+t\zeta)\) with \(\mathbf 1^T\zeta=0\), write
\(G(t)=I+tH+\tfrac12t^2G''(0)+o(t^2)\). The identity \(\operatorname{tr}G\equiv M\) gives
\(\operatorname{tr}G''(0)=0\), hence
\begin{equation}
\label{eq:square-logdet-hessian}
\left.\frac{d^2}{dt^2}\log\det G(t)\right|_{t=0}
=
\operatorname{tr}G''(0)-\operatorname{tr}(H^2)
=
-\|H\|_F^2 .
\end{equation}
With \(H=X\epsilon+Y\zeta\), the right-hand side contains the mixed term
\(-2\langle X\epsilon,Y\zeta\rangle_F\), which must vanish by (i). This is
\eqref{eq:square-tangent-orthogonality}.
\end{remark}

\subsection{Classical square-DFT energy}
For the one-dimensional square DFT, the deficit law reaches its critical-sampling endpoint, and
the irreducible defect is a classical energy.

\begin{proposition}[Square-DFT energy and continuum limit]
\label{prop:square-dft-energy}
Let \(m=s\) in the one-dimensional setting of Section~\ref{sec:one-dim}. For every real node
tangent \(\epsilon\), the fixed-weight and optimally reweighted first-order defects coincide and
equal
\begin{equation}
\label{eq:square-dft-fourier-energy}
q(\epsilon)
=
s\sum_{r=1}^{s-1}
r(s-r)\,|\widehat\epsilon_r|^2
=
E_s(\epsilon),
\end{equation}
with the node-space representation
\begin{equation}
\label{eq:square-dft-energy-identity}
q(\epsilon)
=
\frac14
\sum_{j\neq k}
\frac{|\epsilon_j-\epsilon_k|^2}
{\sin^2\!\bigl(\pi(j-k)/s\bigr)} .
\end{equation}
If \(\epsilon(x)=\sum_{n\in\mathbb Z}\epsilon_n e^{-inx}\in H^\sigma(\mathbb T)\) with
\(\sigma>\tfrac12\) is sampled on the equispaced grid, then
\begin{equation}
\label{eq:square-dft-continuum}
\lim_{s\to\infty}\frac{1}{s^2}q(\epsilon)
=
\sum_{n\neq0}|n|\,|\epsilon_n|^2
=
\|\epsilon\|_{\dot H^{1/2}(\mathbb T)}^2 .
\end{equation}
\end{proposition}

\begin{proof}
At \(m=s\) one has \(\delta=s\) and \(r_j=j\), so Theorem~\ref{thm:deficit} gives
\eqref{eq:square-dft-fourier-energy}; Proposition~\ref{prop:exact-square-rigidity} shows that
the fixed-weight and optimally reweighted defects agree. For
\eqref{eq:square-dft-energy-identity}, expand
\(\epsilon_k=\sum_r\widehat\epsilon_re^{-irx_k}\). For fixed \(d=j-k\),
\(\sum_k|\epsilon_{k+d}-\epsilon_k|^2=4s\sum_r|\widehat\epsilon_r|^2\sin^2(\pi rd/s)\), and the
finite Fej\'er identity
\(\sum_{d=1}^{s-1}\sin^2(\pi rd/s)/\sin^2(\pi d/s)=r(s-r)\) gives the claim.

For the limit, write the normalized form as
\begin{equation}
\label{eq:square-dft-normalized-energy}
\frac{1}{s^2}q(\epsilon)
=
\sum_{r=1}^{s-1}
\frac{r(s-r)}{s}\,
|\widehat\epsilon_r|^2 .
\end{equation}
If \(\epsilon\) is a trigonometric polynomial and \(s\) exceeds twice its degree, each
\(\widehat\epsilon_r\) equals \(\epsilon_n\) for the unique \(n\equiv r\) in \((-s/2,s/2]\), and
\(r(s-r)/s=|n|(1-|n|/s)\), so
\(s^{-2}q(\epsilon)=\sum_{n\ne0}|n|(1-|n|/s)|\epsilon_n|^2\to\sum_{n\ne0}|n||\epsilon_n|^2\).
For general \(\epsilon\in H^\sigma\), the sampled coefficient is the alias sum
\(\widehat\epsilon_r=\sum_{k\in\mathbb Z}\epsilon_{\tilde r+ks}\) with
\(\tilde r\equiv r\) in \((-s/2,s/2]\), and \(r(s-r)/s\le|\tilde r|\). Weighted
Cauchy--Schwarz with weights \(\langle\tilde r+ks\rangle^{\pm2\sigma}\) gives
\[
|\tilde r|\,|\widehat\epsilon_r|^2
\le
|\tilde r|\Bigl(\sum_k\langle\tilde r+ks\rangle^{-2\sigma}\Bigr)
\sum_k\langle\tilde r+ks\rangle^{2\sigma}|\epsilon_{\tilde r+ks}|^2 ,
\]
and the prefactor is bounded by a constant \(C_\sigma\) independent of \(s\): the \(k=0\) term
is at most \(\langle\tilde r\rangle^{1-2\sigma}\le1\), and for \(k\ne0\) one has
\(|\tilde r+ks|\ge|k|s/2\) and \(|\tilde r|\le s/2\), so the remaining terms contribute
\(O(s^{1-2\sigma})\). Summing over \(r\) gives the uniform bound
\[
s^{-2}q(\epsilon)\le C_\sigma\|\epsilon\|_{H^\sigma}^2 .
\]
The limiting form is also bounded by \(\|\epsilon\|_{H^\sigma}^2\). Both are quadratic forms,
so approximating \(\epsilon\) in \(H^\sigma\) by trigonometric polynomials and applying the
triangle inequality to their square roots extends the limit to all of \(H^\sigma\).
\end{proof}

The finite multiplier remains \(r(s-r)/s\); its limit \(|n|\) for fixed signed modes is the
classical Douglas, or circle Dirichlet-to-Neumann, energy
\cite{Douglas1931,DiNezzaPalatucciValdinoci2012}. The node-space identity
\eqref{eq:square-dft-energy-identity} identifies the finite form with the Hessian, at
equispaced nodes, of the circular log-Vandermonde energy
\[
-2\sum_{j<k}\log\left|e^{ix_j}-e^{ix_k}\right| ,
\]
whose \(r(s-r)\) spectrum is classical \cite{CalogeroPerelomov1978,Forrester2010}. What is new
is the identification of these classical energies with the irreducible Fourier-Gram defect.

Since every coefficient \(r(s-r)\) is positive for \(1\le r\le s-1\),
\begin{equation}
\label{eq:square-dft-nullspace}
q(\epsilon)=0
\qquad\Longleftrightarrow\qquad
\widehat\epsilon_r=0
\quad (r=1,\ldots,s-1).
\end{equation}
The only first-order harmless motion at the square DFT is therefore a rigid translation of all
nodes. Critical sampling is the opposite endpoint from the difference-set-injective regime of
Section~\ref{sec:no-collisions}: there, every sufficiently small node motion can be absorbed by
reweighting; here, no non-rigid node motion can be improved by reweighting even to first
order.

\section{What the block geometry looks like in higher dimensions}
\label{sec:higher-dim}

\subsection{Tensor boxes}
\label{sec:tensor-boxes}
For Cartesian frequency boxes on product grids, the one-dimensional alias geometry extends
coordinatewise.
\begin{equation}
\label{eq:tensor-box-geometry}
\Lambda
=
\prod_{a=1}^d\{0,\ldots,m_a-1\},
\qquad
\mathcal G
=
\prod_{a=1}^d\mathbb Z_{s_a},
\qquad
m_a\le s_a .
\end{equation}

For coordinate \(a\), write \(\Delta^{(a)}_{r_a}\) for the corresponding one-dimensional
residue class and
\begin{equation}
\label{eq:tensor-coordinate-multiplicity}
N^{(a)}_{p_a}=m_a-|p_a|
\end{equation}
for its raw-difference multiplicity.

To write the factorization explicitly, let
\begin{equation}
\label{eq:tensor-coordinate-mass}
Z^{(a)}_{r_a}
=
\sum_{p_a\in\Delta^{(a)}_{r_a}}N^{(a)}_{p_a},
\end{equation}
and let \(C^{(a)}_{r_a}\) denote the corresponding one-dimensional covariance in coordinate
\(a\). If \(\Delta^{(a)}_{r_a}\) is empty, we set
\(Z^{(a)}_{r_a}=C^{(a)}_{r_a}=0\). The class data then factor coordinatewise, with each
diagonal entry of \(C_r\) measuring alias spread in one coordinate, weighted by the
multiplicities in the others.

\begin{theorem}[Tensor-box factorization]
\label{thm:tensor-box}
Let \(\Lambda\) and \(\mathcal G\) be the Cartesian frequency box and sampling grid of
\eqref{eq:tensor-box-geometry} with \(m_a\le s_a\), let \(Z^{(a)}_{r_a}\) be as in
\eqref{eq:tensor-coordinate-mass}, and let \(\bar p^{(a)}_{r_a}\) and \(C^{(a)}_{r_a}\) be
the corresponding one-dimensional class center and covariance, all set to zero on empty
classes. Then the class multiplicities factor as
\begin{equation}
\label{eq:tensor-class-multiplicity}
N_r
=
\prod_{a=1}^d Z^{(a)}_{r_a},
\end{equation}
the class centers are coordinatewise,
\[
\bar p_r=\bigl(\bar p^{(1)}_{r_1},\ldots,\bar p^{(d)}_{r_d}\bigr)
\quad\text{when }N_r>0,
\qquad
\bar p_r=0
\quad\text{when }N_r=0,
\]
and every alias-covariance block is diagonal and given by
\begin{equation}
\label{eq:tensor-block-factorization}
C_r
=
\operatorname{diag}_{a=1,\ldots,d}
\left(
C^{(a)}_{r_a}
\prod_{b\neq a}Z^{(b)}_{r_b}
\right).
\end{equation}
Writing
\begin{equation}
\label{eq:tensor-diagonal-entry}
\kappa_a(r)
=
C^{(a)}_{r_a}\prod_{b\neq a}Z^{(b)}_{r_b}
\;\ge\;0,
\end{equation}
the reduced defect separates by coordinate,
\begin{equation}
\label{eq:tensor-defect-coordinatewise}
q(\epsilon)
=
\sum_r\sum_{a=1}^d
\kappa_a(r)\,\bigl|\widehat\epsilon_{r,a}\bigr|^2 ,
\end{equation}
and the harmless polarizations in class \(r\) are
\begin{equation}
\label{eq:tensor-kernel-coordinatewise}
\ker C_r
=
\bigl\{z\in\mathbb C^d: z_a=0\ \text{whenever }\kappa_a(r)>0\bigr\}.
\end{equation}
Moreover \(\kappa_a(r)>0\) precisely when coordinate \(a\) collides in its class, that is
\(|\Delta^{(a)}_{r_a}|\ge2\), and every other coordinate class \(\Delta^{(b)}_{r_b}\) is
nonempty.
\end{theorem}

\begin{proof}
Because both the frequency box and the grid are products, the difference set and its
multiplicities factor,
\[
\Delta=\prod_a\Delta^{(a)},
\qquad
\Delta_r=\prod_a\Delta^{(a)}_{r_a},
\qquad
N_p=\prod_aN^{(a)}_{p_a},
\]
the middle identity because congruence modulo the grid annihilator holds coordinatewise.
Summing \(N_p\) over \(\Delta_r\) therefore factors into one-dimensional sums, which is
\eqref{eq:tensor-class-multiplicity}; if some \(\Delta^{(a)}_{r_a}\) is empty then
\(Z^{(a)}_{r_a}=0\) and \(N_r=0\), as the convention requires. When \(N_r>0\), the first
moment in coordinate \(a\) factors as
\(\bigl(\sum_{p_a}N^{(a)}_{p_a}p_a\bigr)\prod_{b\neq a}Z^{(b)}_{r_b}\), and dividing by
\(N_r\) leaves \(\bar p^{(a)}_{r_a}\); the center is therefore coordinatewise.

The \((a,b)\) entry of \(C_r\) is
\(\sum_{p\in\Delta_r}N_p(p_a-\bar p^{(a)}_{r_a})(p_b-\bar p^{(b)}_{r_b})\). For \(a\neq b\)
the product structure separates this into
\[
\Bigl(\sum_{p_a}N^{(a)}_{p_a}\bigl(p_a-\bar p^{(a)}_{r_a}\bigr)\Bigr)
\Bigl(\sum_{p_b}N^{(b)}_{p_b}\bigl(p_b-\bar p^{(b)}_{r_b}\bigr)\Bigr)
\prod_{c\neq a,b}Z^{(c)}_{r_c},
\]
and each of the first two factors is a centered first moment, hence zero; so \(C_r\) is
diagonal. For \(a=b\) the same separation leaves the one-dimensional centered second moment
\(C^{(a)}_{r_a}\) multiplied by the zeroth moments \(\prod_{b\neq a}Z^{(b)}_{r_b}\), which is
\eqref{eq:tensor-block-factorization} and \eqref{eq:tensor-diagonal-entry}. If \(N_r=0\),
some \(Z^{(a)}_{r_a}\) vanishes, so every diagonal entry vanishes and \(C_r=0\).

Since \(C_r=\operatorname{diag}\bigl(\kappa_1(r),\ldots,\kappa_d(r)\bigr)\) with
\(\kappa_a(r)\ge0\), we have
\(\widehat\epsilon_r^{\,*}C_r\widehat\epsilon_r
=\sum_a\kappa_a(r)|\widehat\epsilon_{r,a}|^2\), and summing over \(r\) in
\eqref{eq:alias-block-decomposition} gives \eqref{eq:tensor-defect-coordinatewise}. A
diagonal positive-semidefinite matrix annihilates exactly those vectors whose components
vanish at the strictly positive diagonal entries, which is
\eqref{eq:tensor-kernel-coordinatewise}. Finally \(\kappa_a(r)>0\) requires
\(C^{(a)}_{r_a}>0\) together with \(Z^{(b)}_{r_b}>0\) for every \(b\neq a\); the latter says
each other coordinate class is nonempty, and the former says \(\Delta^{(a)}_{r_a}\) carries
at least two distinct representatives, because a one-dimensional scatter with positive
multiplicities vanishes exactly on singletons.
\end{proof}

\subsection{Rank-1 lattice}
\label{sec:rank1-lattice}
The block geometry is not peculiar to Cartesian product grids. Consider the two-dimensional
hyperbolic cross
\begin{equation}
\label{eq:rank1-hyperbolic-cross}
\Lambda_R
=
\left\{
k\in\mathbb Z^2:
(1+|k_1|)(1+|k_2|)\le 4
\right\},
\qquad M=17,
\end{equation}
sampled on the rank-1 lattice
\begin{equation}
\label{eq:rank1-lattice}
x_j=\frac{2\pi j}{23}(1,5),
\qquad j=0,\ldots,22 .
\end{equation}
The seventeen resolved frequencies have distinct characters on this lattice, so the sampling
rule is exact on \(\Lambda_R\) \cite{Kaemmerer2013}. Their pairwise differences behave very
differently: the difference set collides in almost every character class. In fact,
\[
|\Delta|=69>23,
\]
so the 69 raw differences must collide among the 23 available characters. Thus
reconstructibility of the resolved frequency set does not imply difference-set injectivity,
which is the stronger property relevant to reweightability.

The collisions are not merely numerous; they can span the full displacement space. In this
example the covariance blocks have rank histogram \(1,2,20\) for ranks \(0,1,2\),
respectively. A representative full-rank class is \(r=2\), with raw differences
\begin{equation}
\label{eq:rank1-full-rank-class}
(-3,1),\qquad (-1,-4),\qquad (0,5),\qquad (2,0),
\end{equation}
having multiplicities \(4,2,2,7\). Their weighted center is
\begin{equation}
\label{eq:rank1-center}
\bar p_2=\left(0,\frac25\right),
\end{equation}
so \(C_2\) is positive definite. A reconstructing rank-1 lattice can still contain alias
classes with no nontrivial first-order harmless polarization.

The nonzero class center \(\bar p_2\) in \eqref{eq:rank1-center} also yields a generally
nonzero optimal weight response through \eqref{eq:optimal-class-weight}. Exact reconstruction
here coexists with both difference collisions and nontrivial first-order reweighting.

Any four points in \(\mathbb R^2\) are affinely dependent, so this class is too, and
Theorem~\ref{thm:exactness-manifold} does not apply. Direct dense evaluation of the reduced defect agrees with the block formula to
relative error below \(10^{-14}\). The example illustrates the covariance-block mechanism
beyond Cartesian grids without asserting a general characterization of rank-1 lattices.

\subsection{Frobenius versus spectral objectives}
\label{sec:frobenius-vs-spectral}
The optimization studied throughout this paper is a Frobenius problem: for a prescribed node
motion, the weights are chosen to minimize
\[
\|X\epsilon+Y\zeta\|_F .
\]
That choice need not minimize worst-case spectral distortion, which we measure for a
Hermitian tangent \(H\) by its spectral half-spread
\(\tfrac12\bigl(\lambda_{\max}(H)-\lambda_{\min}(H)\bigr)\). For example, take the square DFT
with \(s=4\) and
\[
\epsilon=(1,0,-1,0).
\]
By \eqref{eq:square-reweighting-orthogonality}, the Frobenius-optimal correction is
\(\zeta=0\). The resulting tangent matrix \(H\) has characteristic polynomial
\begin{equation}
\label{eq:spectral-counterexample-charpoly}
(\lambda^2-\lambda-1)(\lambda^2+\lambda-1),
\end{equation}
so its spectral half-spread is \((1+\sqrt5)/2\).

A different admissible weight tangent improves the spectral objective. Taking
\[
\zeta=\left(0,\frac14,0,-\frac14\right)
\]
gives
\begin{equation}
\label{eq:spectral-counterexample-square}
(H+Y\zeta)^2=2I,
\end{equation}
so the spectral half-spread drops to
\[
\sqrt2 .
\]
At the same time, however, the squared Frobenius norm increases from \(6\) to \(8\). Hence the
Frobenius-optimal correction can fail to be spectral-spread-optimal. The covariance calculus
therefore answers a specific question: optimal diagonal reweighting for the first-order
Frobenius frame-operator defect. It should not be read as a solution to worst-case spectral
optimization.

\section{Discussion}
\label{sec:discussion}
Alias geometry determines the first-order limits of diagonal density compensation: the
covariance blocks \(C_r\) measure irreducible damage and identify its harmless polarizations.
The tensor DFT shows why this description is not the whole local geometry. Affine relations
among collided differences can obstruct exact deformation even along directions in \(\ker q\).

Several questions remain outside the present theory. When an exact node set admits more than
one positive normalized weight vector, the values and kernels of the reduced form depend on
which exact reference weight is chosen. Section~\ref{sec:no-collisions} shows that this
dependence is an explicit reparametrization, so optimizing over the exact-weight polytope is
tractable as a finite convex quadratic program; what remains open is how the
reference-optimized geometry interacts with the nonlinear exactness set. Likewise, the
analysis begins at an exact rule \(G=I\), so perturbations of already inexact references
require a different expansion. Finally, affine dependence inside collision classes is
understood here only through the second-order compatibility mechanism and the tensor-DFT
example; a general description of the nonlinear exactness set in such classes remains open. First-order alias geometry is the object of study here, and
Proposition~\ref{prop:affine-relation-obstruction} is offered as the entry point for that
higher-order analysis, not as a partial attempt at it.

The covariance viewpoint may nevertheless be useful beyond the examples treated here. It
converts a global first-order reweighting problem into independent low-dimensional blocks
whose ranks, kernels, and eigenvectors identify which perturbation polarizations are harmless
and which are intrinsically damaging. Because those blocks depend only on the geometry and
multiplicities of colliding frequency differences, they can be computed before solving any
particular reconstruction problem and can serve as diagnostics for comparing
sampling geometries, perturbation directions, or candidate lattices. The tensor and rank-1
examples show that this interpretation is not tied to one-dimensional bands or Cartesian
grids; what remains to be understood is how far the same classwise geometry extends once the
exact subgroup structure itself is relaxed.

\begingroup
\sloppy

\endgroup


\begin{thebibliography}{99}
\setlength{\itemsep}{0pt}
\setlength{\parskip}{0pt}

\bibitem{AustinTrefethen2017}
A. P. Austin and L. N. Trefethen,
``Trigonometric interpolation and quadrature in perturbed points,''
\emph{SIAM J. Numer. Anal.}, 55(5):2113--2122, 2017.

\bibitem{Banica2014}
T. Banica,
``First order deformations of the Fourier matrix,''
\emph{J. Math. Phys.}, 55:013504, 2014.

\bibitem{BarrosSaBengtsson2013}
N. Barros e S\'a and I. Bengtsson,
``Families of complex Hadamard matrices,''
\emph{Linear Algebra Appl.}, 438(7):2929--2957, 2013.

\bibitem{BenedettoFickus2003}
J. J. Benedetto and M. Fickus,
``Finite normalized tight frames,''
\emph{Adv. Comput. Math.}, 18(2--4):357--385, 2003.

\bibitem{CalogeroPerelomov1978}
F. Calogero and A. M. Perelomov,
``Properties of certain matrices related to the equilibrium configuration of the one-dimensional many-body problems with the pair potentials \(V_1(x)=-\log|\sin x|\) and \(V_2(x)=1/\sin^2x\),''
\emph{Commun. Math. Phys.}, 59:109--116, 1978.

\bibitem{ChenKutyniokOkoudjouPhilippWang2015}
X. Chen, G. Kutyniok, K. A. Okoudjou, F. Philipp, and R. Wang,
``Measures of scalability,''
\emph{IEEE Trans. Inf. Theory}, 61(8):4410--4423, 2015.

\bibitem{DiNezzaPalatucciValdinoci2012}
E. Di Nezza, G. Palatucci, and E. Valdinoci,
``Hitchhiker's guide to the fractional Sobolev spaces,''
\emph{Bull. Sci. Math.}, 136(5):521--573, 2012.

\bibitem{Dita2004}
P. Di\c{t}\u{a},
``Some results on the parametrization of complex Hadamard matrices,''
\emph{J. Phys. A: Math. Gen.}, 37(20):5355--5374, 2004.

\bibitem{Douglas1931}
J. Douglas,
``Solution of the problem of Plateau,''
\emph{Trans. Amer. Math. Soc.}, 33(1):263--321, 1931.

\bibitem{FeichtingerGroechenigStrohmer1995}
H. G. Feichtinger, K. Gr\"ochenig, and T. Strohmer,
``Efficient numerical methods in non-uniform sampling theory,''
\emph{Numer. Math.}, 69(4):423--440, 1995.

\bibitem{Forrester2010}
P. J. Forrester,
\emph{Log-Gases and Random Matrices},
London Mathematical Society Monographs, Princeton University Press, 2010.

\bibitem{GolubPereyra1973}
G. H. Golub and V. Pereyra,
``The differentiation of pseudo-inverses and nonlinear least squares problems whose variables separate,''
\emph{SIAM J. Numer. Anal.}, 10(2):413--432, 1973.

\bibitem{Groechenig1993}
K. Gr\"ochenig,
``A discrete theory of irregular sampling,''
\emph{Linear Algebra Appl.}, 193:129--150, 1993.

\bibitem{HogeKwanPike1997}
R. D. Hoge, R. K. S. Kwan, and G. B. Pike,
``Density compensation functions for spiral MRI,''
\emph{Magn. Reson. Med.}, 38(1):117--128, 1997.

\bibitem{JonesNjastadThron1989}
W. B. Jones, O. Nj{\aa}stad, and W. J. Thron,
``Moment theory, orthogonal polynomials, quadrature, and continued fractions associated with the unit circle,''
\emph{Bull. Lond. Math. Soc.}, 21(2):113--152, 1989.

\bibitem{Kadec1964}
M. I. Kadets,
``The exact value of the Paley--Wiener constant,''
\emph{Dokl. Akad. Nauk SSSR}, 155(6):1253--1254, 1964;
English transl., \emph{Soviet Math. Dokl.}, 5:559--561, 1964.

\bibitem{Kaemmerer2013}
L. K\"ammerer,
``Reconstructing hyperbolic cross trigonometric polynomials by sampling along rank-1 lattices,''
\emph{SIAM J. Numer. Anal.}, 51(5):2773--2796, 2013.

\bibitem{KieferWolfowitz1960}
J. Kiefer and J. Wolfowitz,
``The equivalence of two extremum problems,''
\emph{Canad. J. Math.}, 12:363--366, 1960.

\bibitem{KircheisPottsFrontiers2023}
M. Kircheis and D. Potts,
``Fast and direct inversion methods for the multivariate nonequispaced fast Fourier transform,''
\emph{Front. Appl. Math. Stat.}, 9:1155484, 2023.

\bibitem{KircheisPottsDCF2023}
M. Kircheis and D. Potts,
``Optimal density compensation factors for the reconstruction of the Fourier transform of bandlimited functions,''
in \emph{Proceedings of the 14th International Conference on Sampling Theory and Applications (SampTA)}, 2023; arXiv:2304.00094.

\bibitem{KutyniokOkoudjouPhilippTuley2013}
G. Kutyniok, K. A. Okoudjou, F. Philipp, and E. K. Tuley,
``Scalable frames,''
\emph{Linear Algebra Appl.}, 438(5):2225--2238, 2013.

\bibitem{Peherstorfer2011}
F. Peherstorfer,
``Positive trigonometric quadrature formulas and quadrature on the unit circle,''
\emph{Math. Comp.}, 80(275):1685--1701, 2011.

\bibitem{PipeMenon1999}
J. G. Pipe and P. Menon,
``Sampling density compensation in MRI: rationale and an iterative numerical solution,''
\emph{Magn. Reson. Med.}, 41(1):179--186, 1999.

\bibitem{Rosenfeld1998}
D. Rosenfeld,
``An optimal and efficient new gridding algorithm using singular value decomposition,''
\emph{Magn. Reson. Med.}, 40(1):14--23, 1998.

\bibitem{SedaratNishimura2000}
H. Sedarat and D. G. Nishimura,
``On the optimality of the gridding reconstruction algorithm,''
\emph{IEEE Trans. Med. Imaging}, 19(4):306--317, 2000.

\bibitem{Simon2005OPUC1}
B. Simon,
\emph{Orthogonal Polynomials on the Unit Circle. Part 1: Classical Theory},
American Mathematical Society Colloquium Publications, vol. 54, Amer. Math. Soc.,
Providence, RI, 2005.

\bibitem{TadejZyczkowski2006}
W. Tadej and K. \.{Z}yczkowski,
``A concise guide to complex Hadamard matrices,''
\emph{Open Syst. Inf. Dyn.}, 13(2):133--177, 2006.

\bibitem{TadejZyczkowski2008}
W. Tadej and K. \.{Z}yczkowski,
``Defect of a unitary matrix,''
\emph{Linear Algebra Appl.}, 429(2--3):447--481, 2008.

\bibitem{YuTownsend2023}
A. Yu and A. Townsend,
``On the stability of unevenly spaced samples for interpolation and quadrature,''
\emph{BIT Numer. Math.}, 63:23, 2023.

\end{thebibliography}
\end{document}